\documentstyle
    [times,amsmath,amssymb,amsthm,pstricks,pst-node,cite,aps,prl,twocolumn,floats]
    {revtex}

\newtheorem{theorem}{Theorem}
\newtheorem{lemma}[theorem]{Lemma}
\newtheorem{question}{Question}

\theoremstyle{remark}
\newtheorem*{exercise}{Exercise}

\newcommand{\R}{{\mathbb{R}}}
\newcommand{\C}{{\mathbb{C}}}
\newcommand{\Z}{{\mathbb{Z}}}
\newcommand{\F}{{\mathbb{F}}}
\newcommand{\Q}{{\mathbb{Q}}}

\newcommand{\tr}{{\mathrm{tr}}}
\newcommand{\Pf}{{\mathrm{Pf}}}
\renewcommand{\sl}{{\mathrm{sl}}}
\DeclareMathOperator{\coker}{coker}
\DeclareMathOperator{\im}{im}

\renewcommand{\tensor}{\otimes}
\newcommand{\ie}{{\em i.e.}}

\newcommand{\vi}{\vec{\imath}}
\newcommand{\vj}{\vec{\jmath}}
\newcommand{\vk}{\vec{k}}
\renewcommand{\tilde}{\widetilde}
\renewcommand{\bar}{\overline}

\newenvironment{fullfigure}[2]
    {\begin{figure}[htb]\begin{center}\def\fullfiga{#1}\def\fullfigb{#2}}
    {\vspace{\baselineskip}\caption{\fullfigb.}\label{\fullfiga}\end{center}\end{figure}}

\newcommand{\Fig}[1]{Figure~\ref{#1}}
\newcommand{\Thm}[1]{Theorem~\ref{#1}}
\newcommand{\Sec}[1]{Section~\ref{#1}}

\psset{linewidth=.5pt,dash=3pt 3pt,doublesep=.05,dotsize=1pt 5}
\SpecialCoor
\newcommand{\psgoesto}{\hspace{.5cm}\pspicture[.5](0,-.1)(1,.1)
    \psline{->}(0,0)(1,0)
    \endpspicture\hspace{.5cm}}

\newcommand{\psvec}[2]{\rput{#1}(#2){\psline{->}(-.15,0)(.15,0)}}
\newcommand{\whitedot}{\pscircle[fillstyle=solid,fillcolor=white](0,0){2.5pt}}
\newcommand{\blackdot}{\qdisk(0,0){2.5pt}}
\newcommand{\whitedisk}{\pscircle[fillstyle=solid,fillcolor=white]}
\newcommand{\blackdisk}{\qdisk}
\newcommand{\daline}{\psline[linestyle=dashed]}
\newcommand{\thline}{\psline[linewidth=1.5pt]}
\newcommand{\dtline}{\psline[linestyle=dashed,linewidth=1.5pt]}
\newcommand{\mto}{\lput{:U}{\pspicture(0,0)(0,0)
\psline[arrows=->,arrowscale=1.5](3.5pt,0)(3.6pt,0)\endpspicture}}
\newcommand{\mfro}{\lput{:U}{\pspicture(0,0)(0,0)
\psline[arrows=->,arrowscale=1.5](-3.5pt,0)(-3.6pt,0)\endpspicture}}

\begin{document}

\wideabs{
\title{An exploration of the permanent-determinant method}
\author{Greg Kuperberg}
\address{UC Davis \\ E-mail: \tt greg@math.ucdavis.edu}
\maketitle
\begin{abstract}
\setlength{\parskip}{.5ex}
The permanent-determinant method and its generalization, the Hafnian-Pfaffian
method, are methods to enumerate perfect matchings of plane graphs that was
discovered by P. W. Kasteleyn.  We present several new techniques and arguments
related to the permanent-determinant with consequences in enumerative
combinatorics.  Here are some of the results that follow from these techniques:

1. If a bipartite graph on the sphere with $4n$ vertices is invariant under the
antipodal map, the number of matchings is the square of the number of matchings
of the quotient graph.

2. The number of matchings of the edge graph of a graph with vertices of degree
at most 3 is a power of 2.

3. The three Carlitz matrices whose determinants count $a \times b \times c$
plane partitions all have the same cokernel.

4. Two symmetry classes of plane partitions can be enumerated with almost no
calculation.
\end{abstract}}

The permanent-determinant method and its generalization, the Hafnian-Pfaffian
method, is a method to enumerate perfect matchings of plane graphs that was
discovered by P. W. Kasteleyn \cite{Kasteleyn:crystal}.  Given a bipartite
plane graph $Z$, the method produces a matrix whose determinant is the number
of perfect matchings of $Z$.  Given a non-bipartite plane graph $Z$, it
produces a Pfaffian with the same property.  The method can be used to
enumerate symmetry classes of plane partitions \cite{Kuperberg:perdet,Kuperberg:oneroof} and domino
tilings of an Aztec diamond \cite{Yang:thesis} and is related to some recent
factorizations of the number of matchings of plane graphs with symmetry
\cite{Jockusch:perfect,Ciucu:reflective}.  It is related to the Gesset-Viennot
lattice path method \cite{GV:partitions}, which has also been used to enumerate
plane partitions \cite{Stembridge:enumeration,Andrews:tsscpp}.  The method
could lead to a fully unified enumeration of all ten symmetry classes of plane
partitions.  It may also lead to a proof of the conjectured $q$-enumeration of
totally symmetric plane partitions.

In this paper, we will discuss some basic properties of the
permanent-determinant method and some simple arguments that use it.  Here are
some original results that follow from the analysis:

\begin{enumerate}
\item If a bipartite graph on the sphere with $4n$ vertices is
invariant under the antipodal map, the number of matchings is the
square of the number of matchings of the quotient graph.
\item The number of matchings of the edge graph of a graph
with vertices of degree at most 3 is a power of 2.
\item The three Carlitz matrices whose determinants count $a \times
b \times c$ plane partitions all have the same cokernel.
\item Two symmetry classes of plane partitions
can be enumerated with almost no calculation.  (This result
was independently found by Ciucu \cite{Ciucu:reflective}).
\end{enumerate}

The paper is largely written in the style of an expository, emphasizing
techniques for using the permanent-determinant method rather than specific
theorems that can be proved with the techniques. Here is a summary for the
reader interested in comparing with previously known results:  Sections
\ref{s:graphs}, \ref{s:perdet}, and \ref{s:hafpfaff} are a review of well-known
linear algebra and results of Kasteleyn, except for \ref{s:spin} and
\ref{s:projective}, which are new. Sections \ref{s:symmetry}, \ref{s:gv}, and
\ref{s:other} are mostly new.  Parts of Section \ref{s:symmetry} were
discovered independently by Regge and Rasetti, Jockusch, Ciucu, and Tesler.
Obviously the Gessel-Viennot method, the Ising model, and tensor calculus
themselves are due to others.  Section \ref{s:pp} consists entirely of new and
independently discovered results about plane partitions. Finally Section
\ref{s:historical} is strictly a historical survey.

\subsection{Acknowledgements}

The author would like to thank Mihai Ciucu and especially Jim Propp for
engaging discussions and meticulous proofreading.  The author also had
interesting discussions about the present work with John Stembridge and Glenn
Tesler.

The figures for this paper were drafted with PSTricks \cite{PSTricks}.
The paper is typeset in two-column REV\TeX with the Times font
package.

\section{Graphs and determinants}
\label{s:graphs}

A {\em sign-ordering} of a finite set is a linear ordering
chosen up to an even permutation.  Given two disjoint sets $A$ and $B$, a
bijection $f:A \to B$ induces a sign-ordering of $A \cup B$ as follows.  Order
the elements of $A$ arbitrarily, and then list
$$a_1,f(a_1),a_2,f(a_2),\ldots.$$
More generally, an oriented matching of a finite set $A$, meaning a partition
of $A$ into ordered pairs, induces a sign-ordering of $A$ by the same
construction. A sign-ordering of $A \cup B$ is also equivalent to a linear
ordering of $A$ and a linear ordering of $B$, chosen up to simultaneous odd or
even permutations, by choosing $f$ to be order-preserving.

Let $Z$ be a weighted bipartite graph with black and white vertices, where the
weights of the edges lie in some field $\F$.  (Usually $\F$ will will be $\R$
or $\C$.)  The graph $Z$ has a {\em weighted, bipartite adjacency matrix},
$M(Z)$, whose rows are indexed by the black vertices of $Z$ and whose columns
are indexed by the white vertices.  The matrix entry $M(Z)_{v,w}$ is the total
weight of all edges from $v$ to $w$.  If the vertices of $Z$ are sign-ordered,
then $\det(M(Z))$ is well-defined (and taken to be 0 unless $M(Z)$ is
square).   By abuse of notation, we define
$$\det(Z) = \det(M(Z)).$$
The sign of $\det(Z)$ is determined by choosing linear orderings of the rows
and columns compatible with the sign-ordering of $Z$.   If the vertices are not
sign-ordered, the absolute determinant $|\det(Z)|$ is still well-defined.

Just as matrices are a notation for linear transformations, a weighted
bipartite graph $Z$ can also denote a linear transformation
$$L(Z):\F[B] \to \F[W].$$
Here $B$ is the set of black vertices, $W$ is the set of white vertices, and
$\F[X]$ denotes the set of formal linear combinations of elements of $X$ with
coefficients in $\F$. The map $L(Z)$ is the one whose matrix is $M(Z)$. Note
that $Z$ is not uniquely determined by $L(Z)$: if $Z$ has multiple edges, the
linear transformation only depends on the sum of the weights of these edges. If
$Z$ has an edge with weight 0, the edge is synonymous with an absent edge.  Row
and column operations on $M(Z)$ can be viewed as operations on $Z$ itself
modulo these ambiguities.

These observations also hold for weighted, oriented non-bipartite graphs. Given
such a graph $Z$, the {\em antisymmetric adjacency matrix} $A(Z)$ has a row and
column for every vertex of $Z$.  The matrix entry $A(Z)_{v,w}$ is the total
weight of all edges from $v$ to $w$ minus the total weight of edges from $w$ to
$v$.  This matrix has a Pfaffian $\Pf(A(Z))$ whose sign is well-defined if the
vertices of $Z$ are sign-ordered.  We also define
$$\Pf(Z) = \Pf(A(Z)).$$

Recall that the Pfaffian $\Pf(M)$ of an antisymmetric matrix $M$
is a sum over matchings in the set of rows of $M$.  The sign of the Pfaffian
depends on a sign-ordering of the rows of $M$.  In these respects,
the Pfaffian generalizes the determinant.  The Pfaffian
also satisfies the relation
\begin{equation}
    \det(M) = \Pf(M)^2. \label{e:detpf}
\end{equation}
This relation has a bijective proof: If $M$
is antisymmetric, the terms in the determinant indexed by permutations with
odd-length cycles vanish or cancel in pairs. The remaining terms are bijective
with pairs of matchings of the rows of $M$, and the signs agree.  This
argument, and the permanent-determinant method generally, blur the distinction
between bijective and algebraic proofs in enumerative combinatorics.

In particular,
    $$\det(Z) = \Pf(Z)$$
when $Z$ is bipartite if all edges are oriented from black to white. (If this
seems inconsistent with equation~\eqref{e:detpf}, recall that the implicit
matrix on the right, $A(Z)$, has two copies of the one on the left, $M(Z)$.) If
$Z$ has indeterminate weights, the polynomial $\det(Z)$ or $\Pf(Z)$ has one
term for each perfect matching $m$ of $Z$.  The term may be written
    $$t(m) = (-1)^m \omega(m),$$
where $(-1)^m$ is the sign of $m$ relative to the sign-ordering of vertices of
$Z$, and $\omega(m)$ is the product of the weights of edges of $m$. Thus
$\det(Z)$ or $\Pf(Z)$ for an arbitrary graph $Z$ is a basic object in
enumerative combinatorics.

\section{The permanent-determinant method}
\label{s:perdet}

Let $Z$ be a connected, bipartite planegraph.  By planarity we mean that $Z$
is embedded in an (oriented) sphere.  The faces of $Z$ are disks; together with
the edges and vertices they form a cell structure, or CW complex, on the
sphere.  Since the sphere is oriented, each face is oriented.  The edges of $Z$
have a preferred orientation, namely the one in which all edges point from
black to white.  The Kasteleyn curvature (curvature for short) of $Z$ at a face
$F$ is defined as
$$c(F) = (-1)^{|F|/2+1} \frac{\prod_{e \in F_+}\omega(e)}
{\prod_{e \in F_-} \omega(e)},$$
where $|F|$ is the number of edges in $F$, $F_+$ is the set of edges whose
orientation agrees with the orientation of $F$, and
$F_-$ is the set of edges whose orientation
disagrees with that of $F$.  (Each face inherits its orientation from
that of the sphere.) A face $F$ is flat if $c(F)$ is 1.  See \Fig{f:weights}.

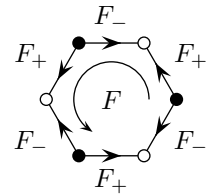
\begin{fullfigure}{f:weights}{Computing Kasteleyn curvature}
\psset{xunit=.433cm,yunit=.75cm}
\pspicture(-3,-1.5)(3,1.5)
\pcline(-1,-1)(-2, 0)\Aput{$F_-$}\mto
\pcline(-2, 0)(-1, 1)\Aput{$F_+$}\mfro
\pcline(-1, 1)( 1, 1)\Aput{$F_-$}\mto
\pcline( 1, 1)( 2, 0)\Aput{$F_+$}\mfro
\pcline( 2, 0)( 1,-1)\Aput{$F_-$}\mto
\pcline( 1,-1)(-1,-1)\Aput{$F_+$}\mfro
\qdisk(-1,-1){2.5pt}\pscircle[fillstyle=solid,fillcolor=white](-2,0){2.5pt}
\qdisk(-1, 1){2.5pt}\pscircle[fillstyle=solid,fillcolor=white]( 1,1){2.5pt}
\qdisk( 2, 0){2.5pt}\pscircle[fillstyle=solid,fillcolor=white](1,-1){2.5pt}
\psarc{->}(0,0){.5}{0}{240}
\rput(0,0){$F$}
\endpspicture
\end{fullfigure}

\begin{theorem}[Kasteleyn] If $Z$ is unweighted, a flat weighting exists.
\label{th:kexists}
\end{theorem}

The theorem depends on the following lemma.

\begin{lemma} If $Z$ has an even number of vertices, and in particular
if black and white vertices are equinumerous, then there are an
even number of faces with $4k$ sides.
\end{lemma}

\begin{proof} Let $n_V$, $n_E$, and $n_F$ be the number of vertices,
edges, and faces of $Z$, respectively.  The Euler characteristic
equation of the sphere is
$$\chi = n_F - n_E + n_V = 2.$$
The term $n_V$ is even.  Divide the contribution to $n_E$ from each edge,
namely $-1$, evenly between the two incident faces. Then the contribution to
$n_F - n_E$ of a face with $4k$ sides is an odd integer, while the contribution
of a face with $4k+2$ sides is an even integer.  Therefore there are an even
number of the former.
\end{proof}

\begin{proof}[Proof of theorem]
Consider the cohomological chain complex of the cell structure given by $Z$
with coefficients in the multiplicative group $\F^*$.  (Since it may
be confusing to consider homological algebra with multiplicative coefficients,
we will sometimes denote a ``sum'' of $\F^*$-cochains as $a\dotplus b$.)
Consider the same orientations of the edges and
faces of $Z$ as above.  With these orientations, we can view a function from
$n$-cells to $\F^*$ as an $n$-cochain. In particular, a weighting $\omega$ of
$Z$ is equivalent to a $1$-cochain.  Let $\omega_k$, the Kasteleyn cochain, be
a 2-cochain which assigns $(-1)^{|F|/2+1}$ to each face $f$.  The coboundary
$\delta \omega$ of $\omega$ is related to the curvature by
$$c(F) = \omega_k \dotplus \delta \omega.$$
Thus, a flat weighting exists if and only if $\omega_k$ is a coboundary.  By
the lemma, $\omega_k$ has an even number of faces with weight $-1$ and the rest
have weight 1.  Thus, $\omega_k$ represents the trivial second cohomology class
of the sphere.  Therefore it is a coboundary.
\end{proof}

Following the terminology of the proof, the curvature of any weighting is a
coboundary, because it is the sum (in the sense of ``$\dotplus$'') of two
coboundaries, $\omega_k$ and the coboundary of the weighting.  Thus the product
of all curvatures of all faces is 1.

\begin{theorem}[Kasteleyn] If $Z$ is flat, the number
of perfect matchings is $\pm\det(Z)$, because $t(m)$ has
the same sign for all $m$. \label{th:perdet}
\end{theorem}

A complete proof is given in Reference \citen{Kuperberg:perdet}, but the result
also follows from a more general result.

By a {\em loop} we mean a collection of edges of $Z$ whose union is a simple
closed curve. If the loop $\ell$ is the difference between two matchings $m_1$
and $m_2$, then all edges of $\ell$ point in the same direction if we reverse
the edges of $\ell \cap m_2$.  Of the two regions of the sphere separated by
$\ell$, the positive one is the one whose orientation agrees with $\ell$.

\begin{theorem} If $m_1$ and $m_2$ are two matchings of $Z$ that
differ by one loop $\ell$, the ratio of their terms $t(m_1)/t(m_2)$
in the expansion of $\det(Z)$ equals the product of the 
curvatures of the faces on the positive side of $\ell$.  \label{th:curvdet}
\end{theorem}

\begin{proof} The loop $\ell$ has an even number of sides and also must
enclose an even number vertices on the positive side $S_+$.  If we remove the
vertices and edges on the negative side $S_-$, we obtain a new graph $Z'$
such that the loop $\ell$ bounds a face $F$ that replaces $S_-$. Since the
total curvature of all faces of $Z'$ is 1, the curvature of $F$
is the reciprocal of the total curvature of all other faces.  Finally,
\begin{align*}
c(F)\ & =(-1)^{|F|/2+1}
\frac{\prod_{e \in F_+} \omega(e)}{\prod_{e \in F_-} \omega(e)} \\
& = (-1)^{m_1} (-1)^{m_2} \frac{\prod_{e \in \ell \cap m_2} \omega(e)}
{\prod_{e \in \ell \cap m_1} \omega(e)} \\
& = \frac{t(m_2)}{t(m_1)}.
\end{align*}
The signs agree because $m_1$ and $m_2$ differ by an even cycle, which is an
odd permutation if and only if $F$ has $4k$ sides.
\end{proof}

\begin{fullfigure}{f:loop}{A loop enclosing four faces}
\pspicture(-1.5,-1.5)(1.5,1.5)
\pnode(1.5; 22.5){A1} \pnode(1.5; 67.5){A2} \pnode(1.5;112.5){A3}
\pnode(1.5;157.5){A4} \pnode(1.5;202.5){A5} \pnode(1.5;247.5){A6}
\pnode(1.5;292.5){A7} \pnode(1.5;337.5){A8}
\pnode(0,.6){B1}\pnode(-.45,-.375){B2}

\psline(A1)(B1)(A3)
\psline(A4)(B2)(A6)
\thline(B1)(B2)
\pspolygon[linewidth=1.5pt](A1)(A2)(A3)(A4)(A5)(A6)(A7)(A8)
\rput(1.05;67.5){$F_1$}
\rput(.9;135){$F_2$}
\rput(1.05;202.5){$F_3$}
\rput(.525;337.5){$F_4$}
\rput(-1.7,0){$\ell$}
\qdisk(A1){2.5pt} \qdisk(A2){2.5pt} \qdisk(A3){2.5pt} \qdisk(A4){2.5pt}
\qdisk(A5){2.5pt} \qdisk(A6){2.5pt} \qdisk(A7){2.5pt} \qdisk(A8){2.5pt}
\qdisk(B1){2.5pt} \qdisk(B2){2.5pt}

\endpspicture
\end{fullfigure}

\Fig{f:loop} illustrates the proof of \Thm{th:curvdet}.
The loop encloses four faces.  Edges in bold appear in at least one of
two terms $m_1$ and $m_2$ that differ by $\ell$.  The theorem
in this case says that
$$\frac{t(m_2)}{t(m_1)} = c(F_1)c(F_2)c(F_3)c(F_4).$$

In light of \Thm{th:curvdet}, if $Z$ is an unweighted graph, a curvature
function and a reference matching $m$ are enough to define $\det(Z)$, because
we can choose a weighting and a sign-ordering with the desired curvature and
such that $t(m) = 1$.  The matrix $M(Z)$ will then have the following
ambiguity.  In general, if $\omega_1$ and $\omega_2$ are two weightings with
the same curvature, then $\omega_1\dotplus\omega_2^{-1}$ is a 1-cocycle.  Since
the first homology of the sphere is trivial, the ratio is a 1-coboundary, \ie,
$\omega_1$ and $\omega_2$ differ by a 0-cochain.  The corresponding matrices
$M(Z_{\omega_1})$ and $M(Z_{\omega_2})$ then differ by multiplication by
diagonal matrices on the left and the right.

\section{The Hafnian-Pfaffian method}
\label{s:hafpfaff}

Kasteleyn's method for non-bipartite plane graphs expresses the number of
perfect matchings as a Pfaffian.  For simplicity, we consider unweighted,
oriented graphs. The analysis has a natural generalization to weighted graphs
in which the orientation is completely separate from the weighting. The {\em
curvature} of an orientation at a face is 1 if an odd number of edges point
clockwise around the face and $-1$ otherwise.  The orientation is flat if the
curvature is 1 everywhere. A routine generalization of \Thm{th:perdet}
shows that if $Z$ has a flat orientation, $\Pf(Z)$ is the number of matchings
\cite{Kasteleyn:crystal}.

A graph $Z$, whether planar or not, is {\em Pfaffian} if it admits an
orientation such that all terms in $\Pf(Z)$ have the same sign
\cite{Little:k33free}.

If $Z$ is bipartite, orienting $Z$ is equivalent to giving each edge
a 1 if it points to black to white and $-1$ otherwise. The weighting is flat if
and only if the orientation is flat.

\begin{theorem}[Kasteleyn] If $Z$ is an unoriented plane graph, a
flat orientation exists. \label{th:kpfexists} 
\end{theorem}

In particular, planar graphs are Pfaffian.

\begin{proof} The proof follows that of \Thm{th:kexists}. Fix an
orientation $o$, and again consider the mod 2 Euler characteristic of the
sphere.  Ignoring the vertices, we transfer the Euler characteristic of each
edge to the incident face whose orientation agrees with that of the edge. 
The net Euler characteristic of a face is then 0 if it is flat and 1 if it is
not, therefore there must be an even number of non-flat faces.  Let $k$ be
the curvature of $o$.

The Euler characteristic calculation shows that $k$ is a coboundary of a
1-cochain $c$ with coefficients in the multiplicative group $\{\pm 1\}$.  Let
$o' = c\dotplus o$ be the ``sum'' of $c$ and $o$, defined by the rule that $o'$ and $o$
agree on those edges where $c$ is $1$ and disagree where $c$ is $-1$.  Then
$o'$ is flat.
\end{proof}

In the same vein, suppose that $Z$ and $Z'$ are the same graph with two
different flat orientations.  By homology considerations, the ``difference'' of
the two orientations (1 where they agree, $-1$ where they disagree) is a
1-coboundary.  Thus they differ by the coboundary of a 0-cochain $c$, which is
a function from the vertices to $\{\pm 1\}$.  Let $D$ be the diagonal matrix
whose entries are the values of $c$.  Then $A(Z)$, $A(Z')$, and $D$ satisfy the
relation
$$A(Z) = DA(Z')D^T.$$
Note that $c$ and $D$ are unique up to sign.

\subsection{Spin structures}
\label{s:spin}

We conclude with some comments about flat orientations of a graph $Z$ on a
surface of genus $g$.  Kasteleyn \cite{Kasteleyn:crystal} proved that the
number of matchings of such a graph is given by a sum of $4^g$ Pfaffians
defined using inequivalent flat orientations of $Z$.  (See also Tesler
\cite{Tesler:non-oriented}.) There is an interesting relationship between these
flat orientations and spin structures.  A spin structure on a surface is
determined by a vector field with even-index singularities. We can make such a
vector field using an orientation and a matching $m$.  At each vertex, make the
vectors point to the vertex.  Then replace each edge by a continuous family of
edges such that in the middle of the edge, the vector field is 90 degrees
clockwise relative to the orientation of the edge \Fig{f:vectors} shows this
operation applied to the four edges of a square.

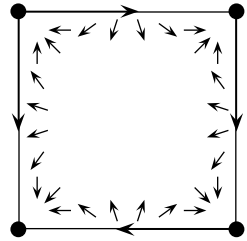
\begin{fullfigure}{f:vectors}{Vectors describing a spin structure}
\pspicture(-1.6,-1.45)(1.6,1.45)
\psframe(-1.45,-1.45)(1.45,1.45)
\psline[arrowsize=3pt 4]{->}(-1.45,1.45)(  .15, 1.45)
\psline[arrowsize=3pt 4]{->}(1.45, 1.45)( 1.45, -.15)
\psline[arrowsize=3pt 4]{->}(1.45,-1.45)( -.15,-1.45)
\psline[arrowsize=3pt 4]{->}(-1.45,1.45)(-1.45, -.15)
\qdisk(-1.45,-1.45){3pt}
\qdisk( 1.45,-1.45){3pt}
\qdisk( 1.45, 1.45){3pt}
\qdisk(-1.45, 1.45){3pt}
\psvec{ 90}{ 1.2,  .9} \psvec{ 45}{ 1, 1} \psvec{  0}{  .9, 1.2}
\psvec{180}{ -.9, 1.2} \psvec{135}{-1, 1} \psvec{ 90}{-1.2,  .9}
\psvec{270}{-1.2, -.9} \psvec{225}{-1,-1} \psvec{180}{ -.9,-1.2}
\psvec{  0}{  .9,-1.2} \psvec{315}{ 1,-1} \psvec{270}{ 1.2, -.9}

\psvec{324}{ .54,1.2} \psvec{288}{ .18,1.2}
\psvec{252}{-.18,1.2} \psvec{216}{-.54,1.2}

\psvec{126}{-1.2, .54} \psvec{162}{-1.2, .18}
\psvec{198}{-1.2,-.18} \psvec{234}{-1.2,-.54}

\psvec{ 36}{ .54,-1.2} \psvec{ 72}{ .18,-1.2}
\psvec{108}{-.18,-1.2} \psvec{144}{-.54,-1.2}

\psvec{126}{1.2, .54} \psvec{162}{1.2, .18}
\psvec{198}{1.2,-.18} \psvec{234}{1.2,-.54}

\endpspicture
\end{fullfigure}

Because the orientation is flat, the vector field extends to the faces with
even-index singularities, but the singularities at the vertices are odd. 
Contract the odd-index singularities in pairs along edges of the matching; the
resulting vector field induces a spin structure.  For a fixed orientation,
inequivalent matchings yield distinct spin structures.  Here two matchings are
equivalent if they are homologous.  For a fixed matching, two inequivalent
orientations yield distinct spin structures.

\subsection{Projective-plane graphs}
\label{s:projective}

An expression for the number of matchings of a non-planar graph may in general
require many Pfaffians.  But there is an interesting near-planar case when a
single Pfaffian suffices.

A graph is a {\em projective-plane graph} if it is embedded in the projective plane.
A graph embedded in a surface is {\em locally bipartite} if all faces are disks
and have an even number of sides.  It is {\em globally bipartite} if it is
bipartite.  If $Z$ is locally but not globally bipartite, then it has a
non-contractible loop, but all non-contractible loops have odd length while all
contractible loops have even length.

\begin{theorem} If $Z$ is a connected, projective-plane graph
which is locally but not
globally bipartite, then it is Pfaffian.
\end{theorem}

\begin{proof}
Assume that $Z$ has an even number of vertices.

The curvature of an orientation of $Z$ is well-defined even though
the projective plane is non-orientable:  Since each face has an even number of
sides, the curvature is 1 if an odd number of edges point in both directions
and $-1$ otherwise. If the curvature of an arbitrary orientation $o$ is a
coboundary, meaning that an even number of faces have curvature $-1$, then there
is a flat orientation by the homology argument of \Thm{th:kpfexists}.

To prove that the curvature of $o$ must be a coboundary, we cut along a
non-contractible loop $\ell$, which must have odd length, to obtain an
oriented plane graph $Z'$. Every face of $Z$ becomes a face of $Z'$, and in addition
$Z'$ has an outside face with $2|\ell|$ sides.  A face in $Z'$ has the same
curvature as in $Z$ assuming that it is a face of both graphs. The graph $Z'$
has an odd number of vertices, because it has $|\ell|$ more vertices than
$Z$ does. By the argument of \Thm{th:kpfexists}, $Z'$ has an odd number
of faces with curvature $-1$.  Moreover, the outside face is one of them,
because each of the edges of $\ell$ appears twice, both times pointing either
clockwise or counterclockwise.  Therefore $Z$ must have an even number of faces
with curvature $-1$.

Finally, we show that a flat orientation of $Z$ is in fact Pfaffian. Let $m_1$
and $m_2$ be two matchings that differ by a single loop. Since the loop has
even length, it is contractible.  By the usual argument, the ratio
$t(m_1)/t(m_2)$ of the
corresponding terms in the expansion of $\Pf(A(Z))$ equals the product of the
curvatures of the faces that the loop bounds.  Since $Z$
is flat, this product is 1.
\end{proof}

\section{Symmetry}
\label{s:symmetry}

\subsection{Generalities}
\label{s:generalities}

Let $V$ be a vector space over $\C$, the complex numbers. If a linear
transformation $L:V \to V$ commutes with the action of a reductive group $G$,
then $\det(L)$ factors according to the direct sum decomposition of $V$ into
irreducible representations of $G$.  At the abstract level, for each distinct
irreducible representation $R$, we can make a vector space $V_R$ such that $V_R
\tensor R$ is an isotypic summand of $V$, and there exist isotypic blocks
$$L_R:V_R \to V_R$$
such that
$$L = \bigoplus_R (L_R \tensor I_R),$$
where $I_R$ is the identity on $R$. Then
\begin{equation}
\det(L) = \prod_R \det(L_R)^{\dim R}. \label{e:factor}
\end{equation}
More concretely, if $M$ is a matrix and a group $G$ has a matrix representation
$\rho$ such that
$$\rho(g)M = M \rho(g),$$
then after a change of basis, $M$ decomposes into blocks, with $\dim R$
identical blocks of some size for each irreducible representation $R$ of $G$,
so its determinant factors.

Suppose that $L$ is an endomorphism of some integral lattice $X$ in $V$
(concretely, if $M$ is an integer matrix) and $R$ is some rational
representation.  After choosing a rational basis $\{r_i\}$ for $R$, we can
realize copies $V_{r_i}$ of $V_R$ as rational subspaces of $V$.  The lattice
$L$ preserves each $V_{r_i}$ and acts acts on it as $L_R$.  Then $X \cap
V_{r_i}$ is a lattice in $V_{r_i}$, and $L$ is an endomorphism of this lattice
as well.  The conclusion is that each $\det(L_R)$ must be an integer because
$L_R$ is an endomorphism of a lattice. Indeed, this argument works for any
number field (such as the Gaussian rationals) and its ring of integers (such as
the Gaussian integers) if $R$ is not a rational representation, which tells us
that equation~\ref{e:factor} is in general a factorization into algebraic
integers if $L$ is integral.  The determinant $\det(L_R)$ is, a priori, in the
same field as the representation $R$.   A refinement of the argument shows that
it is in the same field as the character of $R$, which may lie in a smaller
field than $R$ itself.

\begin{fullfigure}{f:pfaffian}{A Pfaffian orientation with broken symmetry}
\pspicture(-1.2,-1.2)(1.2,1.2)
\pcline( .5, .5)(-.5, .5)\mfro
\pcline(-.5, .5)(-.5,-.5)\mto
\pcline(-.5,-.5)( .5,-.5)\mfro
\pcline( .5,-.5)( .5, .5)\mfro

\pcline( 1.2, 1.2)(-1.2, 1.2)\mfro
\pcline(-1.2, 1.2)(-1.2,-1.2)\mto
\pcline(-1.2,-1.2)( 1.2,-1.2)\mfro
\pcline( 1.2,-1.2)( 1.2, 1.2)\mfro

\pcline(-.5,-.5)(-1.2,-1.2)\mfro
\pcline( .5,-.5)( 1.2,-1.2)\mto
\pcline( .5, .5)( 1.2, 1.2)\mfro
\pcline(-.5, .5)(-1.2, 1.2)\mto
\endpspicture
\end{fullfigure}

The general principle of factorization of determinants applies to enumeration
of matchings in graphs with symmetry via the Hafnian-Pfaffian method. As
discussed in Sections~\ref{s:graphs} and \ref{s:hafpfaff}, an oriented graph
$Z$ yields an antisymmetric map
$$A(Z):\C[Z] \to \C[Z].$$
Any symmetry of the oriented graph intertwines this map, and the factorization
principle applies. However there are three complications.  First, the
orientation may have less symmetry than the graph itself (\Fig{f:pfaffian}).
Second, the principle of factorization gives information about the determinant
and not the Pfaffian.  (If a summand $R$ of an orthogonal representation $V$ is
orthogonal and $L$ is an antisymmetric endomorphism of $V$, then the factor
$\det(L_R)$ is the square of $\Pf(L_R)$.  But if $R$ is symplectic or complex,
then $\det(L_R)$ need not be a square.  In these cases the factorization
principle is less informative.) Third, the number of matchings might only
factor into algebraic integers, which is less informative than a factorization
into ordinary integers. 

Let $Z$ be a connected plane graph, and suppose that a group $G$ acts on the
sphere and preserves $Z$ and the orientation of the sphere. Then $G$ acts by
permutation matrices on $V = \C[Z]$, the vector space generated by vertices of
$Z$.  Although $G$ commutes with the adjacency matrix of $Z$, it does not in
general commute with the antisymmetric adjacency matrix $A(Z)$ if $Z$ is
oriented, because $G$ might not preserve the orientation of $Z$.  However, if
$Z$ has a flat orientation $o$ and $g \in G$, then $go$ differs from $o$ by a
coboundary in the sense of Sections~\ref{s:perdet} and \ref{s:hafpfaff}. This
means that there is a signed permutation matrix $\tilde{g}$ which does commute
with $A(Z)$.  These signed permutation matrices together form a linear
representation of some group $\tilde{G}$ which extends $G$. At first glance it
may appear as if the fiber of this extension is as big as $\{\pm 1\}^|Z|$.  But
because $Z$ is connected, among diagonal signed matrices only  the identity and
its negative commute with $A(Z)$; only constant 0-cochains leave alone the
orientation of every edge of $Z$. Therefore $\tilde{G}$ is a central extension
of $G$ by the two-element group $G_0 = \{\pm 1\}$. The subgroup $G_0$ either
acts trivially on $V$ or negates it.  Thus, in decomposing $V$ into
irreducible representations, we need only consider those where $G_0$ acts
trivially (by definition the {\em even} representations) or only those where
$G_0$ acts by negation (by definition the {\em odd} representations), depending
on whether the action on $V$ is even or odd.

If $G$ has odd order, the central extension must split. In this case, by
averaging, $Z$ has a flat orientation invariant under $G$
\cite{Tesler:private}.  If $G$ a cyclic group of rotations of order $2n$, then
the central extension might not split, but it is not very interesting as an
extension; if $Z$ is bipartite, one can find a flat weighting using $4n$th
roots of unity which is invariant under $G$ \cite{Jockusch:perfect}.  But in
the most complicated case, when $Z$ has icosahedral symmetry, the central
extension $\tilde{G}$ (when it is non-trivial) is the binary icosahedral group
$\tilde{A_5}$, which is quite interesting.  This central extension seems
related to the connection between flat orientations and spin structures
mentioned in Section~\ref{s:hafpfaff}, because the symmetry group of a spun
sphere is an analogous central extension of $SO(3)$. Irrespective of $G$, the
representation theory of $\tilde{G}$ reveals a factorization of the number of
matchings of $Z$.

If $Z$ is bipartite, then there are two important changes to the story. First,
after including signs, one can make orientation-reversing symmetries commute
with $A(Z)$ as well, because in the bipartite case they take flat orientations
to flat orientations. If $Z$ is not bipartite, the best signed versions of
orientation-reversing symmetries instead anticommute with $A(Z)$. 
Anticommutation is less informative than commutation, but they still sometimes
provide information together with the following fact from linear algebra:  If
$A$ and $B$ anticommute and $B$ is invertible, then
$$\tr(A^n) = \tr(-A^n) = 0$$
for $n$ odd, because
$$-A^n = B^{-1}A^nB.$$

Second, if $Z$ is bipartite, then color-reversing symmetries yield no direct
information via representation theory.  In this case, it is better to apply
representation theory to the bipartite adjacency matrix $M(Z)$.  This
matrix represents a linear map
$$L:\C[B] \to \C[W],$$
where $B$ is the set of black vertices and $W$ is the set of white vertices,
rather than a linear endomorphism of a single space. The color-preserving
symmetries act on both $V = \C[B]$ and $U = \C[W]$ and $M(Z)$ intertwines these
actions.  Hence for each irreducible $R$ there is an isotypic block
$$L_R:V_R \to U_R.$$
Hence for each $R$ we must choose volume elements on $V_R$ and $U_R$ so that
$\det L_R$ is well-defined. Nevertheless, equation \ref{e:factor} still holds
if it is properly interpreted.

\subsection{Cyclic symmetry}
\label{s:cyclic}

Suppose that $G$ is generated by a single rotation $g$ of order $n$, and let
$\omega$ be an $n$th root of unity when $n$ is odd or $G$ fixes a vertex (the
split case), or an odd $2n$th root of unity when $n$ is even and $G$ does not
fix a vertex (the non-split case).  Then the vertex space $\C[Z]$ has an
isotypic summand $\C[Z]_\omega$.  A suitable set of vectors of the form
$$v+\omega gv + \ldots + \omega^n g^nv,$$
where $v$ is a vertex of $Z$, form a basis of $\C[Z]_\omega$ (assuming certain
sign conventions in the non-split case).  In this case the isotypic blocks of
$A(Z)$ are all represented by weighted plane graphs whose Kasteleyn curvature
can be easily derived from that of $Z$ \cite{Jockusch:perfect}.

If $Z$ is bipartite, then a reflection symmetry produces a similar
factorization, and again the resulting matrices are represented by plane
graphs \cite{Ciucu:reflective}.

The other possibility for a color-preserving cyclic symmetry is a glide-reflection.   In particular, $Z$ may be invariant under
the antipodal map on the sphere.  In this case, the blocks of $M(Z)$ cannot be
represented by plane graphs.  Instead, they produce projective-plane graphs.  So
the number of matchings factors, but the permanent-determinant
method does not identify either factor as an unweighted enumeration.

\subsection{Color-reversing symmetry}
\label{s:reverse}

If $Z$ is bipartite, then a color-reversing symmetry does not a priori lead to
an interesting factorization of the number of matchings of $Z$.  For example,
if the symmetry is an involution, then if we use it to establish a bijection
between black and white vertices, we learn only that $M(Z)$ is symmetric, which
says little about its determinant.  However, color-reversing symmetries do have
two interesting and related consequences.

First, if $Z$ has color-reversing and color-preserving symmetries, the
color-reversing symmetries sometimes imply that the numerical factors of
$\det(Z)$ from the color-preserving symmetries lie in smaller-than-expected
number-field rings.  For example, if $Z$ has a color-reversing 90-degree
rotational symmetry $g$, then the symmetry $g^2$ yields
$$\det(Z) = \det(Z_i)\det(Z_{-i}),$$
where $Z_{\pm i}$ is the quotient graph $Z/g^2$ with curvature $\pm i$ at the
faces fixed by $g^2$.  (Recall that $Z$ is on the sphere, so there are two
fixed faces.)  Then the remaining symmetry tells us that  $Z_i$ and $Z_{-i}$
have a curvature-preserving isomorphism, which means that their determinants
are equal up to a unit in the Gaussian integers.  At the same time, their
determinants are complex conjugates.  Thus, $\det(Z_i)$ is, up to a unit, in
the form $a$ or $(1+i)a$ for some rational integer $a$. The conclusion is that
$\det(Z)$ is either a square or twice a square \cite{Jockusch:perfect}. 
Similarly, if $Z$ has a color-reversing 60-degree symmetry,
$$\det(Z) = ab^2$$
for some integers $a$ and $b$ coming from enumerations in quotient graphs.

Second, if $Z$ has a color-reversing involution $g$ which does not fix any
edges, then the antisymmetric adjacency matrix $A(Z/g)$ of the quotient graph
$Z/g$ can be interpreted as the bipartite adjacency matrix of $Z$ with some
weighting. Since the determinant is the square of the Pfaffian,
$$\det(Z) = \Pf(Z/g).$$
If $Z/g$ is flat (which implies that $Z/g$ has an even number of vertices),
then $Z$ with its induced weighting may or may not be flat, depending on $g$. 
Assuming $Z$ is connected, $g$ may be the antipodal involution in the sphere or
it may be rotation by 180 degrees.  In the first case, $Z$ is flat, while $Z/g$
is projective-plane graph which is locally but not globally bipartite.  In the
second case, $Z$ is not flat, but has curvature $-1$ at the two faces fixed by
rotation.  The first case is a new theorem:

\begin{theorem} If $Z$ is a bipartite graph on the sphere with $4n$ vertices
which is invariant under the antipodal involution $g$, and if $g$ exchanges
colors of vertices of $Z$, then the number of matchings of $Z$ is the square of
the number of matchings of $Z/g$.  \label{th:antipode}
\end{theorem}

For example, the surface of a Rubik's cube satisfies these conditions
(\Fig{f:rubik}).

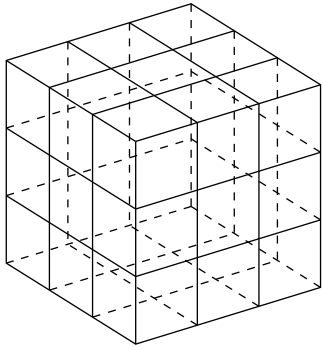
\begin{fullfigure}{f:rubik}{The Rubik's cube graph}
\psset{xunit=.5cm,yunit=.5cm}
\pspicture(-5,-5)(5,5)
\pnode( 0.737, 4.528){A0} \pnode( 0.737,-0.865){A1} \pnode( 4.178, 2.373){A2}
\pnode( 4.178,-3.019){A3} \pnode(-4.178, 3.019){A4} \pnode(-4.178,-2.373){A5}
\pnode(-0.737, 0.865){A6} \pnode(-0.737,-4.528){A7} \pnode(-0.902, 4.025){B0}
\pnode(-0.902,-1.368){B1} \pnode( 2.540, 1.871){B2} \pnode( 2.540,-3.522){B3}
\pnode(-2.540, 3.522){B4} \pnode(-2.540,-1.871){B5} \pnode( 0.902, 1.368){B6}
\pnode( 0.902,-4.025){B7} \pnode( 1.884, 3.810){C0} \pnode( 1.884,-1.583){C1}
\pnode( 3.031, 3.092){C2} \pnode( 3.031,-2.301){C3} \pnode(-3.031, 2.301){C4}
\pnode(-3.031,-3.092){C5} \pnode(-1.884, 1.583){C6} \pnode(-1.884,-3.810){C7}
\pnode( 0.737, 2.730){D0} \pnode( 0.737, 0.933){D1} \pnode( 4.178, 0.576){D2}
\pnode( 4.178,-1.222){D3} \pnode(-4.178, 1.222){D4} \pnode(-4.178,-0.576){D5}
\pnode(-0.737,-0.933){D6} \pnode(-0.737,-2.730){D7}
\daline(A1)(A0) \daline(A1)(A3) \daline(A1)(A5)
\daline(B3)(B1)(B0) \daline(B7)(B5)(B4) \daline(C5)(C1)(C0) \daline(C7)(C3)(C2)
\daline(D2)(D0)(D4) \daline(D3)(D1)(D5)
\psline(A6)(A2) \psline(A6)(A4) \psline(A6)(A7)
\psline(B0)(B2)(B3) \psline(B4)(B6)(B7) \psline(C0)(C4)(C5) \psline(C2)(C6)(C7)
\psline(D4)(D6)(D2) \psline(D5)(D7)(D3)
\pspolygon(A0)(A2)(A3)(A7)(A5)(A4)
\endpspicture
\end{fullfigure}

\begin{exercise} Prove \Thm{th:antipode} with an explicit bijection.
\end{exercise}

This exercise is a special case of the bijective argument
that
$$\det(M) = \Pf(M)^2$$
for any antisymmetric matrix $M$.

\subsection{Icosahedral symmetry}
\label{s:icosahedral}

The results of this section were discovered independently by Rasetti
and Regge \cite{RR:icosahedral}.

The easiest realization of the binary icosahedral group $\tilde{A_5}$ is as the
subgroup of the unit quaternions $a+b\vi + c\vj + d\vk$ for which $(a,b,c,d)$
is one of the points
$$\frac12(\tau,1,\frac1\tau, 0) \qquad (1,0,0,0) \qquad \frac12(1,1,1,1)$$
or the
points obtained from these by changing signs or even permutations of
coordinates. Here $\tau$ is the golden ratio. Note that two elements
of $\tilde{A_5}$ are conjugate if and only if they have the
same real part $a$.

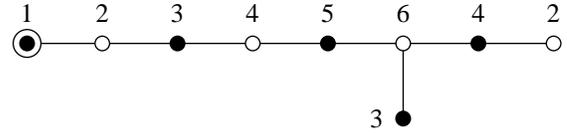
\begin{fullfigure}{f:e8}{Extended $E_8$, a graph of
representations of $\tilde{A_5}$}
\pspicture(0,-.2)(7,1.5)
\psline(0,1)(7,1)
\psline(5,1)(5,0)
\whitedisk(0,1){5.5pt}
\blackdisk(0,1){3pt} \whitedisk(1,1){3pt} \blackdisk(2,1){3pt}
\whitedisk(3,1){3pt} \blackdisk(4,1){3pt} \whitedisk(5,1){3pt}
\blackdisk(6,1){3pt} \whitedisk(7,1){3pt} \blackdisk(5,0){3pt}
\rput(0,1.4){1} \rput(1,1.4){2} \rput(2,1.4){3} \rput(3,1.4){4}
\rput(4,1.4){5} \rput(5,1.4){6} \rput(6,1.4){4} \rput(7,1.4){2}
\rput(4.65,0){3}
\endpspicture
\end{fullfigure}

This realization also describes a two-dimensional representation $\pi$ of
$\tilde{A_5}$.  The character of $\pi$ is twice the real parts of the elements
of $\tilde{A_5}$.  By the McKay correspondence, the irreducible representations
of $\tilde{A_5}$ together form an $E_8$ graph, where the trivial representation
is the extending vertex, and two representations $R$ and $R'$ are joined by
an edge if $R \tensor \pi$ contains $R'$ as a summand.  This
diagram is given in \Fig{f:e8} together with the dimensions of
the representations.  The trivial representation is circled.
A black vertex is an even representation
in the sense of \Sec{s:generalities}, while a white vertex is an odd
representation.
This graph can be used to compute the character table of $\tilde{A_5}$,
which is given in Table~\ref{t:character}.  In this table,
$$\bar{\tau} = -\frac1\tau$$
is the Galois conjugate of $\tau$. The table indicates various properties of
the representations. The conjugacy class $c_0$ contains only $1$, so its row is
the trace of the identity or the dimension of each representation.  The
conjugacy class $c_8$ contains only $-1$, so its row indicates which
representations are even and which are odd. The representation $R_1$ equals the
defining representation $\pi$. Apparently, five of the characters are rational,
while the other four lie in the golden field $\Q(\tau)$. Less superficially,
the character table can be used to find the direct sum decomposition of an
arbitrary representation from its character, or to decompose an equivariant map
into its isotypic blocks.

\begin{table}[htb]
$$\begin{array}{r|rrrrrrrrr}
\chi & R_0 & R_1 & R_2 & R_3 & R_4 & R_5 & R_6 & R_7 & R_8 \\ \hline
c_0 & 1 & 2 & 3 & 4 & 5 & 6 & 4 & 2 & 3 \\
c_1 & 1 & \tau & \tau & 1 & 0 & -1 & -1 & \bar{\tau} & \bar{\tau} \\
c_2 & 1 & 1 & 0 & -1 & -1 & 0 & 1 & 1 & 0 \\
c_3 & 1 & -\bar{\tau} & \bar{\tau} & -1 & 0 & 1 & -1 & -\tau & \tau \\
c_4 & 1 & 0 & -1 & 0 & 1 & 0 & 0 & 0 & -1 \\
c_5 & 1 & \bar{\tau} & \bar{\tau} & 1 & 0 & -1 & -1 & \tau & \tau \\
c_6 & 1 & -1 & 0 &1 & -1 & 0 & 1 & -1 & 0 \\
c_7 & 1 & -\tau & \tau & -1 & 0 & 1 & -1 & -\bar{\tau} & \bar{\tau} \\
c_8 & 1 & -2 & 3 & -4 & 5 & -6 & 4 & -2 & 3
\end{array}$$
\caption{Character table of $\tilde{A_5}$}
\label{t:character}
\end{table}

Suppose that a graph $Z$ has icosahedral symmetry. If a rotation by 180 degrees
fixes a vertex of $Z$, then the action of $\tilde{A_5}$ on $\C[Z]$ is even, but if such a
rotation fixes an edge or a face, then the action is odd (exercise). In the
second case, the factorization principle says that $\det(A(Z))$ is in the form
$a^2\bar{a}^2b^4c^6$, where $b$ and $c$ are integers and $a$ and $\bar{a}$ are
conjugate elements in $\Z(\tau)$, because the available representations have
dimensions 2,2,4, and 6.  Thus the number of matchings factors as
$a^1\bar{a}^1b^2c^3$.  In the first case, $\C[Z]$ decomposes entirely into
orthogonal representations, which implies that
$$\Pf(Z) = a^1b^3\bar{b}^3c^4d^5$$
using the dimensions of the representations and the number fields in which they
lie.   Here are four interesting examples. In the first three, the author
explicitly computed the factorization from symmetry by computing the trace of
$gM(Z)^n$ for different $g$ and $n$.

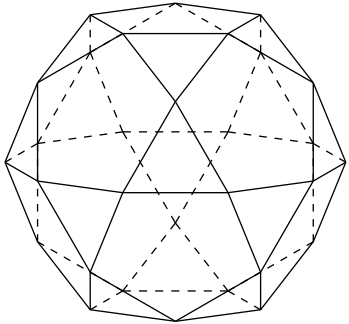
\begin{fullfigure}{f:dodecedge}{The edge graph of a dodecahedron}
\psset{xunit=.7,yunit=.7}
\pspicture(-3.5,-3.5)(3.5,3.5)
\pnode( 0.000, 1.154){Eac} \pnode( 2.618, 1.511){Eae}
\pnode( 2.618,-0.357){Eag} \pnode( 1.000, 2.446){Eai}
\pnode( 1.000,-0.578){Eaj} \pnode( 0.000,-1.154){Ebd}
\pnode( 2.618, 0.357){Ebe} \pnode( 2.618,-1.511){Ebg}
\pnode( 1.000, 0.578){Ebk} \pnode( 1.000,-2.446){Ebl}
\pnode(-2.618, 1.511){Ecf} \pnode(-2.618,-0.357){Ech}
\pnode(-1.000, 2.446){Eci} \pnode(-1.000,-0.578){Ecj}
\pnode(-2.618, 0.357){Edf} \pnode(-2.618,-1.511){Edh}
\pnode(-1.000, 0.578){Edk} \pnode(-1.000,-2.446){Edl}
\pnode( 3.236, 0.000){Eeg} \pnode( 1.618, 2.803){Eei}
\pnode( 1.618, 2.089){Eek} \pnode(-3.236, 0.000){Efh}
\pnode(-1.618, 2.803){Efi} \pnode(-1.618, 2.089){Efk}
\pnode( 1.618,-2.089){Egj} \pnode( 1.618,-2.803){Egl}
\pnode(-1.618,-2.089){Ehj} \pnode(-1.618,-2.803){Ehl}
\pnode( 0.000, 3.023){Eik} \pnode( 0.000,-3.023){Ejl}
\daline(Ebd)(Ebk) \daline(Ebd)(Ebl) \daline(Ebd)(Edk) \daline(Ebd)(Edl)
\daline(Ebe)(Ebg) \daline(Ebe)(Ebk) \daline(Ebe)(Eeg) \daline(Ebe)(Eek)
\daline(Ebg)(Ebl) \daline(Ebk)(Edk) \daline(Ebk)(Eek) \daline(Ebl)(Edl)
\daline(Ebl)(Egl) \daline(Edf)(Edh) \daline(Edf)(Edk) \daline(Edf)(Efh)
\daline(Edf)(Efk) \daline(Edh)(Edl) \daline(Edk)(Efk) \daline(Edl)(Ehl)
\daline(Eei)(Eek) \daline(Eek)(Eik) \daline(Efi)(Efk) \daline(Efk)(Eik)
\psline(Eac)(Eai) \psline(Eac)(Eaj) \psline(Eac)(Eci) \psline(Eac)(Ecj)
\psline(Eae)(Eag) \psline(Eae)(Eai) \psline(Eae)(Eeg) \psline(Eae)(Eei)
\psline(Eag)(Eaj) \psline(Eag)(Eeg) \psline(Eag)(Egj) \psline(Eai)(Eci)
\psline(Eai)(Eei) \psline(Eaj)(Ecj) \psline(Eaj)(Egj) \psline(Ebg)(Eeg)
\psline(Ebg)(Egl) \psline(Ecf)(Ech) \psline(Ecf)(Eci) \psline(Ecf)(Efh)
\psline(Ecf)(Efi) \psline(Ech)(Ecj) \psline(Ech)(Efh) \psline(Ech)(Ehj)
\psline(Eci)(Efi) \psline(Ecj)(Ehj) \psline(Edh)(Efh) \psline(Edh)(Ehl)
\psline(Eei)(Eik) \psline(Efi)(Eik) \psline(Egj)(Egl) \psline(Egj)(Ejl)
\psline(Egl)(Ejl) \psline(Ehj)(Ehl) \psline(Ehj)(Ejl) \psline(Ehl)(Ejl)
\endpspicture
\end{fullfigure}

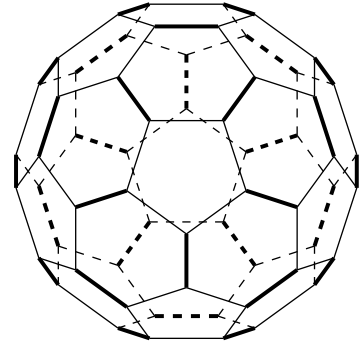
\begin{fullfigure}{f:fullerene}{The $C_{60}$ graph with two kinds of edges}
\pspicture(-2.3,-2.3)(2.3,2.3)
\pnode(-0.420,-1.927){Eac} \pnode( 0.420,-1.927){Eca} \pnode( 1.962,-0.196){Eae}
\pnode( 1.703,-0.995){Eea} \pnode( 1.962,-1.153){Eag} \pnode( 1.703,-1.510){Ega}
\pnode( 0.490,-0.674){Eai} \pnode( 0.910,-1.252){Eia} \pnode( 0.490,-2.223){Eaj}
\pnode( 0.910,-2.086){Eja} \pnode(-0.420, 1.927){Ebd} \pnode( 0.420, 1.927){Edb}
\pnode( 1.962, 1.153){Ebe} \pnode( 1.703, 1.510){Eeb} \pnode( 1.962, 0.196){Ebg}
\pnode( 1.703, 0.995){Egb} \pnode( 0.490, 2.223){Ebk} \pnode( 0.910, 2.086){Ekb}
\pnode( 0.490, 0.674){Ebl} \pnode( 0.910, 1.252){Elb} \pnode(-1.962,-0.196){Ecf}
\pnode(-1.703,-0.995){Efc} \pnode(-1.962,-1.153){Ech} \pnode(-1.703,-1.510){Ehc}
\pnode(-0.490,-0.674){Eci} \pnode(-0.910,-1.252){Eic} \pnode(-0.490,-2.223){Ecj}
\pnode(-0.910,-2.086){Ejc} \pnode(-1.962, 1.153){Edf} \pnode(-1.703, 1.510){Efd}
\pnode(-1.962, 0.196){Edh} \pnode(-1.703, 0.995){Ehd} \pnode(-0.490, 2.223){Edk}
\pnode(-0.910, 2.086){Ekd} \pnode(-0.490, 0.674){Edl} \pnode(-0.910, 1.252){Eld}
\pnode( 2.265,-0.221){Eeg} \pnode( 2.265, 0.221){Ege} \pnode( 0.793, 0.258){Eei}
\pnode( 1.472, 0.478){Eie} \pnode( 0.793, 1.806){Eek} \pnode( 1.472, 1.312){Eke}
\pnode(-2.265,-0.221){Efh} \pnode(-2.265, 0.221){Ehf} \pnode(-0.793, 0.258){Efi}
\pnode(-1.472, 0.478){Eif} \pnode(-0.793, 1.806){Efk} \pnode(-1.472, 1.312){Ekf}
\pnode( 0.793,-1.806){Egj} \pnode( 1.472,-1.312){Ejg} \pnode( 0.793,-0.258){Egl}
\pnode( 1.472,-0.478){Elg} \pnode(-0.793,-1.806){Ehj} \pnode(-1.472,-1.312){Ejh}
\pnode(-0.793,-0.258){Ehl} \pnode(-1.472,-0.478){Elh} \pnode( 0.000, 1.548){Eik}
\pnode( 0.000, 0.834){Eki} \pnode( 0.000,-0.834){Ejl} \pnode( 0.000,-1.548){Elj}
\dtline(Eac)(Eca) \dtline(Eae)(Eea) \dtline(Eai)(Eia) \dtline(Ecf)(Efc)
\dtline(Eci)(Eic) \dtline(Eei)(Eie) \dtline(Eek)(Eke) \dtline(Efi)(Eif)
\dtline(Efk)(Ekf) \dtline(Eik)(Eki) \daline(Eac)(Eic) \daline(Eac)(Ejc)
\daline(Eae)(Ege) \daline(Eae)(Eie) \daline(Eai)(Eci) \daline(Eai)(Eei)
\daline(Ebe)(Eke) \daline(Ebk)(Eek) \daline(Eca)(Eia) \daline(Eca)(Eja)
\daline(Ecf)(Ehf) \daline(Ecf)(Eif) \daline(Eci)(Efi) \daline(Edf)(Ekf)
\daline(Edk)(Efk) \daline(Eea)(Ega) \daline(Eea)(Eia) \daline(Eei)(Eki)
\daline(Eek)(Eik) \daline(Efc)(Ehc) \daline(Efc)(Eic) \daline(Efi)(Eki)
\daline(Efk)(Eik) \daline(Eie)(Eke) \daline(Eif)(Ekf) \thline(Eag)(Ega)
\thline(Eaj)(Eja) \thline(Ebd)(Edb) \thline(Ebe)(Eeb) \thline(Ebg)(Egb)
\thline(Ebk)(Ekb) \thline(Ebl)(Elb) \thline(Ech)(Ehc) \thline(Ecj)(Ejc)
\thline(Edf)(Efd) \thline(Edh)(Ehd) \thline(Edk)(Ekd) \thline(Edl)(Eld)
\thline(Eeg)(Ege) \thline(Efh)(Ehf) \thline(Egj)(Ejg) \thline(Egl)(Elg)
\thline(Ehj)(Ejh) \thline(Ehl)(Elh) \thline(Ejl)(Elj) \psline(Eag)(Eeg)
\psline(Eag)(Ejg) \psline(Eaj)(Ecj) \psline(Eaj)(Egj) \psline(Ebd)(Ekd)
\psline(Ebd)(Eld) \psline(Ebe)(Ege) \psline(Ebg)(Eeg) \psline(Ebg)(Elg)
\psline(Ebk)(Edk) \psline(Ebl)(Edl) \psline(Ebl)(Egl) \psline(Ech)(Efh)
\psline(Ech)(Ejh) \psline(Ecj)(Ehj) \psline(Edb)(Ekb) \psline(Edb)(Elb)
\psline(Edf)(Ehf) \psline(Edh)(Efh) \psline(Edh)(Elh) \psline(Edl)(Ehl)
\psline(Eeb)(Egb) \psline(Eeb)(Ekb) \psline(Efd)(Ehd) \psline(Efd)(Ekd)
\psline(Ega)(Eja) \psline(Egb)(Elb) \psline(Egj)(Elj) \psline(Egl)(Ejl)
\psline(Ehc)(Ejc) \psline(Ehd)(Eld) \psline(Ehj)(Elj) \psline(Ehl)(Ejl)
\psline(Ejg)(Elg) \psline(Ejh)(Elh)
\endpspicture
\end{fullfigure}

\begin{enumerate}
\item If $Z$ is an icosahedron, then $\C[Z]$ consists of  two copies of the
6-dimensional representation $R_6$, which means that $A(Z)$ is, after a change
of basis, six copies of a $2 \times 2$ matrix $M$.  Moreover, $A(Z)$
anticommutes with antipodal inversion, so $M$ must have vanishing trace. The
matrix $M^2$, and therefore $A(Z)$ must be a multiple of the identity. In fact,
$$A(Z)^2 = 5I.$$
Thus there are $5^3 = 125$ matchings.
\item If $Z$ is a dodecahedron, then
$$\Pf(Z) = 36 = 1^3 6^2.$$
The space $\C[Z]$ has two 6-dimensional summands,
which contribute the factors of 1, and two 4-dimensional summands,
which contribute the factors of 6.
\item If $Z$ is the edge graph of a dodecahedron or icosahedron
(\Fig{f:dodecedge}), then
$$\Pf(Z) = (4+2\sqrt5)^3 (4-2\sqrt5)^3 1^4 2^5 = -2^{11}.$$
The space $\C[Z]$ has two of each of the non-trivial even representations.
\item Let $Z$ be the bond graph of the fullerene $C_{60}$ (\Fig{f:fullerene}). 
Suppose that its hexagonal edges (the thicker ones in the figure) have weight 1
and its pentagonal edges (the thinner ones) have weight $p$.  According to
Tesler \cite{Tesler:kekule}, the total weight of all matchings is
\begin{align*}
\Pf(Z) =\ &(1-2p+\frac{5+\sqrt5}2p^2)(1-2p+\frac{5-\sqrt5}2p^2) \nonumber \\
& (1+2p+2p^3+5p^4)^2(1+p^2+2p^3+p^4)^3.
\end{align*}
In this case the space $\C[Z]$ is 60-dimensional and decomposes as two copies
of each 2-dimensional representation, four copies of the odd 4-dimensional
representation, and six copies of the 6-dimensional representation.  The
factors from the 2-dimensional representation are at most quadratic and the
factors of the 4-dimensional representation are at most quartic.  It follows
that the factorization above coincidences with the  factorization given by the
symmetry principle.
\end{enumerate}

\section{Gessel-Viennot}
\label{s:gv}

Some of the results in the section were independently discovered
by Horst Sachs et al \cite{AS:linear}.

The Gessel-Viennot method is another method that counts combinatorial objects
using determinants \cite{GV:partitions}.  Let $Z$ be a directed, plane graph
with $n$ sources (univalent vertices with outdegree 1) and $n$ sinks (univalent
vertices with indegree 1).  Suppose further that the edges at each vertex are
segregated, meaning that no four edges alternate in, out, in, out. Then the
Gessel-Viennot method produces an $n \times n$ matrix whose determinant is the
number of collections of $n$ disjoint, directed paths in $Z$ from the sources
to the sinks.  Each entry of the matrix is the number of paths from a source to
a sink.  The method has some ad hoc generalizations that produce Pfaffians
\cite{Okada:generating,Stembridge:paths}. It has been used to enumerate several
classes of plane partitions \cite{Stembridge:enumeration,Andrews:tsscpp}.

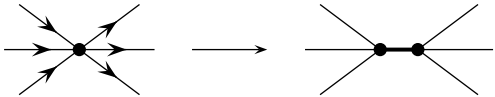
\begin{fullfigure}{f:gvsplit}{Transforming Gessel-Viennot to Kasteleyn}
\pspicture[.5](-1,-.7)(1,.7)
\qdisk(0,0){2.5pt}
\pcline(-.8,-.6)(0,0)\mto \pcline(-.8,.6)(0,0)\mto \pcline( -1,0)(0,0)\mto
\pcline(0,0)(.8,.6)\mto   \pcline(0,0)(1,0)\mto    \pcline(0,0)(.8,-.6)\mto
\endpspicture
\psgoesto
\pspicture[.5](-1,-.7)(1.5,.7)
\psline(-.8,-.6)(0,0)  \psline(-.8, .6)(0,0) \psline(-1,0)(0,0)
\thline(.5,0)(0,0)     \psline(.5,0)(1.3,.6) \psline(.5,0)(1.5,0) 
\psline(.5,0)(1.3,-.6) \qdisk( 0,0){2.5pt}   \qdisk(.5,0){2.5pt}
\endpspicture
\end{fullfigure}

Given a graph $Z$ suitably decorated for the Gessel-Viennot method, there is a
related graph $Z'$ to which the permanent-determinant method applies.  Namely,
split each vertex of $Z$ into an edge $e$, with all inward arrows at one end of
$e$ and all outward arrows at the other end (\Fig{f:gvsplit}). This operation
induces a bijection betwen disjoint path collections in $Z$ and matchings in
$Z'$.

There is a corresponding relation between the Gessel-Viennot matrix $GV(Z)$ and
a Kasteleyn matrix $M(Z')$.  Define a {\em pivot operation} to be the act of
replacing an $n \times n$ matrix of the form
$$\left(\begin{array}{c|c}
M & v \\ \hline w & 1 \end{array}\right)$$
by $M-(v \tensor w)$, where $v$ and $w$ are vectors and $M$ is an $(n-1) \times
(n-1)$ matrix.  The determinant does not change under pivot operations.  Starting
with $M(Z')$, if we pivot at all entries corresponding to edges in $Z'$ which
are contracted in $Z$, the result is the Gessel-Viennot matrix $GV(Z)$ with
some rows and columns negated.  This proves that $\det GV(Z)$ is the number of
matchings of $Z'$.  It also suggests that any enumeration derived using the
Gessel-Viennot method can also be understood using the permanent-determinant
method.

It was pointed out to the author by Mihai Ciucu that there is 
a version of the Gessel-Viennot method for arbitrary graphs, whether
planar or not.  In this case the transformation of Figure~\ref{f:gvsplit}
still produces the relation
$$\det GV(Z) = \det Z',$$
where the right side is a weighted enumeration of matchings of $Z'$ in our
sense.  However, this more general setting cannot be interpreted via
\Thm{th:perdet}.

\subsection{Cokernels}
\label{s:cokernels}

An integer $n \times n$ matrix $M$ can be interpreted as a homomorphism from
$\Z^n$ to itself.  The cokernel is defined as the target divided by the image:
$$\coker M = \Z^n/(\im M).$$
Alternatively, the cokernel is the abelian group on $n$ generators whose
relations are given by $M$.  The cokernel is invariant under pivot operations,
and the number of elements in the cokernel is $|\det(M)|$.  Moreover, the
cokernel of a Kasteleyn matrix $M(Z)$ depends only on the unweighted graph $Z$
and not the particular choice of a flat weighting, and any corresponding
Gessel-Viennot matrix also has the same cokernel.

\begin{question} Is there a natural bijection, or an algebraic generalization
of a bijection, between the cokernel of $M(Z)$ and the set of matching of $Z$?
\end{question}

For any integer matrix $M$, the cokernel of $M^T$ is naturally the Pontryagin
dual of the cokernel of $M$.  In other words, there is a natural Fourier
transform map
$$\Psi:\C[\coker M] \to \C[\coker M^T].$$
This map can be interpreted as a generalized bijection.

\section{Other tricks}
\label{s:other}

\subsection{Forcing planarity}
\label{s:force}

Let $Z$ be an arbitrary bipartite graph with some weighting. Then $\det(Z)$ is
interesting as a weighted enumeration of the matchings of $Z$ and as an
algebraic quantity.  At the same time, there is a simple way to convert $Z$ to
a plane graph without changing its determinant.  Then ideas related to the
permanent-determinant method apply.  (Note that the number of matchings
of $Z$ in general will change.  Together with the assumption
that $NP \ne P$, this would otherwise contradict the fact that
the number of matchings for a general non-planar graph is a $\#P$-hard
quantity \cite{Valiant:permanent}.)

\begin{fullfigure}{f:split}{Vertex splitting}
\pspicture[.5](-.75,-.6)(.75,.6)
\psline(-.6,-.45)(0,0) \psline(-.6,.45)(0,0) \psline( -.75,0)(0,0)
\psline(0,0)(.6,.45)   \psline(0,0)(.75,0)    \psline(0,0)(.6,-.45)
\qdisk(0,0){2.5pt}
\endpspicture
\psgoesto
\pspicture[.5](-.75,-.6)(1.75,.6)
\psline(-.6,-.45)(0,0)  \psline(-.6, .45)(0,0) \psline(-.75,0)(0,0)
\psline(0,0)(1,0)     \psline(1,0)(1.6,.45) \psline(1,0)(1.75,0) 
\psline(1,0)(1.6,-.45)
\qdisk( 0,0){2.5pt}   \qdisk(.5,0){2.5pt}  \qdisk(1,0){2.5pt}
\endpspicture
\end{fullfigure}

In such a graph $Z$, we can {\em triple} an edge, \ie, replace it by three
edges in series.  If the weight of the original edge is $w$ and the weights of
three edges replacing it are $a$, $b$, and $c$, then $\det(Z)$ will not change
if $w = ac$ and $b=-1$.  Edge tripling is a special case of {\em vertex
splitting}, a more general operation in which a vertex is replaced by two edges
in series (\Fig{f:split}).  Vertex splitting changes $\det(Z)$ by a factor
of $\pm a$ if one of the new edges has weight $a$ and the other has weight
$-a$.

\begin{fullfigure}{f:butterfly}{A crossing replaced by a butterfly}
\pspicture(-2,-.75)(2,.75)
\rput(-1.5,0){\pspicture(-.5,-.5)(.5,.5)
\psline(-.5,-.5)(.5, .5) \psline(-.5, .5)(.5,-.5)
\whitedisk(-.5,-.5){2.5pt} \whitedisk(-.5, .5){2.5pt}
\blackdisk( .5, .5){2.5pt} \blackdisk( .5,-.5){2.5pt}
\endpspicture}
\psline{->}(-.5,0)(.5,0)
\rput(5.25,0){\pspicture(-.8,-.75)(8,.75)
\pnode( .7, .6){A1} \pnode( .7,-.6){A2} \pnode(-.7,-.6){A3} \pnode(-.7, .6){A4}
\pnode( .3,0){B1} \pnode(-.3,0){B2}
\pspolygon(B1)(A1)(A4)(B2)(A3)(A2)
\thline(B1)(B2)
\thline(A4)(A1)
\blackdisk(A1){2.5pt} \whitedisk(A4){2.5pt} \blackdisk(B2){2.5pt}
\whitedisk(A3){2.5pt} \blackdisk(A2){2.5pt} \whitedisk(B1){2.5pt}
\endpspicture}
\endpspicture
\end{fullfigure}

Pick a projection of $Z$ in the plane, \ie, a drawing where edges are straight
but may cross.  After tripling the edges sufficiently, every edge crosses at
most one other edge, and the convex hull of two crossing edges contain no
vertices other than the endpoints of those edges.  Then two edges that cross
can be replaced by seven that do not cross according to \Fig{f:butterfly}. Call
the new subgraph a {\em butterfly}.  If the top two horizontal edges
of the butterfly
(which are in bold in the figure) are given weight $-1$ and all other
edges of the butterfly and the edges that cross have weight $1$, then
the operation of replacing the crossing by the
butterfly can be reproduced by row and column operations on $M(Z)$. It does not
change the determinant if the edges are given suitable weights.  Thus, at the
expense of more vertices and edges, $Z$ becomes planar.

\subsection{The Ising model}
\label{s:ising}

The partition function of the unmagnetized Ising model on a graph is the
following weighted enumeration.  Given two numbers $a$ and $b$, called
Boltzmann weights, and given a graph $Z$, compute the total weight of all
functions $s$ (states) from the vertices of $Z$ to the set of  spins
$\{\uparrow,\downarrow\}$, where the weight of a state $s$ is the product of
the weights of the edges, and the weight of an edge is $a$ if the spins of its
vertices agree and $b$ if they disagree.  In a generalization of the model, $a$
and $b$ can be different for different edges.  If $a=b=1$ for a given edge in
this generalization, the edge can be ignored.

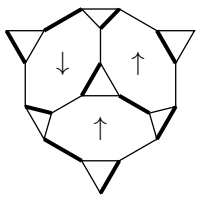
\begin{fullfigure}{f:ising}{An Ising state and its matching}
\psset{xunit=1cm,yunit=.577cm}
\pspicture(-.25,-.5)(2.25,4)
\pnode(   1, -.5){A10}
\pnode( .75, .25){A12} \pnode( .25, .75){A21}
\pnode(1.25, .25){A16} \pnode(1.75, .75){A61}
\pnode(-.25,3.25){A30}
\pnode(   0, 1.5){A23} \pnode(   0, 2.5){A32}
\pnode( .35,1.35){A27} \pnode( .75,1.75){A72}
\pnode( .25,3.25){A34} \pnode( .75,3.75){A43}
\pnode(1.25,3.75){A45} \pnode(1.75,3.25){A54}
\pnode(   1, 3.3){A47} \pnode(   1, 2.5){A74}
\pnode(2.25,3.25){A50}
\pnode(   2, 2.5){A56} \pnode(   2, 1.5){A65}
\pnode(1.65,1.35){A67} \pnode(1.25,1.75){A76}

\rput(1.5,2.5){$\uparrow$}
\rput( .5,2.5){$\downarrow$}
\rput(  1, 1){$\uparrow$}

\pspolygon(A10)(A12)(A16) \pspolygon(A23)(A27)(A21) \pspolygon(A30)(A32)(A34)
\pspolygon(A43)(A45)(A47) \pspolygon(A50)(A54)(A56) \pspolygon(A61)(A65)(A67)
\pspolygon(A72)(A74)(A76)
\thline(A12)(A21) \psline(A16)(A61) \psline(A23)(A32) \psline(A27)(A72)
\thline(A34)(A43) \psline(A45)(A54) \psline(A47)(A74) \psline(A56)(A65)
\thline(A67)(A76)

\thline(A72)(A74)
\thline(A27)(A23)
\thline(A30)(A32)
\thline(A45)(A47)
\thline(A54)(A56)
\thline(A65)(A61)
\thline(A10)(A16)
\endpspicture
\end{fullfigure}

If $Z$ is a plane graph, we can make $Z$ triangular, meaning that all faces are
triangles, at the expense of adding ignorable edges. Let $Z'$ be the dual graph
with a vertex added in the middle of each edge.   Let $Z''$ be the edge graph
of $Z'$.  Then (exercise) there is a 2-to-1 map from the Ising states of $Z$ to
the perfect matchings of $Z''$.  The weights of Ising states can be matched up
to a global factor by assigning weights to edges of $Z''$.  Thus the
Hafnian-Pfaffian method can be used to find the total weight of all Ising
states of $Z$.

More generally, if $Z$ is any plane graph whose vertices have valence 1, 2, or
3, then the number of matchings of the edge graph $Z'$ of $Z$ is a power of
two. Indeed, the set of matchings can always be interpreted as an affine vector
space over $\Z/2$.  A matching of $Z'$ is equivalent to an orientation of $Z$
such that at each vertex, an odd number of edges are oriented outward
(exercise).  If the orientations of each individual edge are arbitrarily
labelled 0 and 1, the constraints at the vertices are all linear. For example,
the number of matchings of the edge graph of a dodecahedron
(Section~\ref{s:icosahedral}) is a power of two.  The bachelorhood vertex
\cite{Kuperberg:perdet} for symmetric plane partitions is a similar
construction.

\subsection{Tensor calculus}
\label{s:tensor}

This section relates determinants of graphs to the formal setting of quantum
link invariants \cite{RT:ribbon}.  We present no complete mathematical results,
only a brief summary of how to discuss determinants of graphs in
more algebraic terms.  Another class of enumeration problems, planar ice and
alternating-sign matrices, are related to the Jones polynomial and other
quantum invariants based on the quantum group $U_q(\sl(2))$.

The Ising model is an example of a {\em state model}, a general scheme where
one has a set of atoms, a set of states for each atom, local weights which
depend on the states of particular clusters of atoms, and the global weight of
a state which is defined as the product of all local weights.  (Most weighted
enumerations of objects such as matchings and plane partitions can be described
as state models.)  Many natural state models can be interpreted as tensor
expressions.  In index notation, one can write an expressions of tensors over a
common vector space $V$ such as
$$A^{ab}_cB^c_dC_b^d\ldots$$
Each index may only appear twice, once covariantly (as a subscript) and once
contravariantly (as a superscript).  Repeated indices are summed.   For
simplicity, suppose all indices are repeated so that the expression is
scalar-valued.  Although this description of index notation depends on choosing
a basis for $V$, such a tensor expression is basis-independent.  In state model
terms, each index is called an atom and the matrix of each tensor is a set of
local weights.  Note also that there an oriented graph associated to a tensor
expression, where each vertex corresponds to a tensor and each edge connecting
to vertices corresponds to indices connecting the tensors.   A tensor
expression is a more versatile form for a state model because arbitrary an
arbitrary linear transformation of $V$ extends to transformations of tensors. 
In particular, $V$ might be a group representation and the tensors might be
invariant under the group action.

If $Z$ is a bipartite, trivially weighted graph, then $\det(Z)$ has such an
interpretation, except that the tensors are invariant under a supergroup action
instead of a group action.  Recall that a supervector space is a
$\Z/2$-graded vector space.  An algebra structure on such a vector
space is supercommutative if it is graded-commutative (odd vectors anticommute
with each other, even-graded vectors commute with everything).  A supergroup is
a group object in the category of supercommutative algebras.  It is like a
group except that instead of a commutative algebra of real- or complex-valued
functions on the group (with multiplication meaning pointwise multiplication),
there is a supercommutative algebra $A$.  The group law (in the
finite-dimensional case) is expressed by a multiplication operation
$$m:A^*\tensor A^* \to A^*$$
on the dual vector space. The particular supergroup of interest to us has the
property that both $A$ and $A^*$ are isomorphic to $\Lambda(\R)$, an exterior
algebra with one generator $x$.  We consider tensors and tensor expressions
over $A$ viewed as a representation of itself.

The space $A$ has an invariant multilinear function
$$y_n:A^{\tensor n} \to \R$$
given by
$$a_1 \tensor a_2 \tensor \ldots \tensor a_n \mapsto \mu(a_1a_2\ldots a_n),$$
where the dual vector
$$\mu:A \to \R$$
is defined by
$$\mu(1) = 0 \qquad \mu(x) = 1.$$
The multilinear function $y_n$ can be viewed as an element of $(A^*)^{\tensor
n}$.  Similarly there is a dual tensor
$$x_n \in A^{\tensor n}$$
using the algebra structure on $A^*$.  If each black vertex of $Z$ is replaced
by a copy of $x_n$ and each white vertex by a copy of $y_n$ in such a way that
there is an index for each edge, then the tensors together form a
scalar-valued expression, because each index appears twice and every index is
contracted.  This scalar-valued expression is exactly $\det(Z)$.

\section{Plane partitions}
\label{s:pp}

Plane partitions are one of the most interesting applications of the
permanent-determinant method \cite{Kuperberg:perdet}. A {\em plane partition in
a box} is a collection of unit cubes in an $a \times b \times c$ box
(rectangular prism) such that below, behind, and to the left of each cube is
either another cube or a wall. A plane partition in a box is equivalent to a
tiling of a hexagon with sides of length $a$, $b$, $c$, $a$, $b$, and $c$ by
unit $60^\circ$ triangles (\Fig{f:pptil}).

\begin{fullfigure}{f:pptil}{A plane partition and a tiling}
\psset{xunit=.5cm,yunit=.577cm}
\pspicture(-7,0)(7,5)
\rput[b](-4.5,0){\pspicture(-2,0)(2,5)
\psset{fillstyle=solid}
\pspolygon[fillcolor=lightgray](-1,4.5)(-1,3.5)(0,3)(0,4)
\pspolygon[fillcolor=darkgray](1,4.5)(1,3.5)(0,3)(0,4)
\pspolygon[fillcolor=lightgray](-1,3.5)(-1,2.5)(0,2)(0,3)
\pspolygon[fillcolor=lightgray](0,3)(0,2)(1,1.5)(1,2.5)
\pspolygon[fillcolor=darkgray](2,3)(2,2)(1,1.5)(1,2.5)
\pspolygon[fillcolor=lightgray](0,2)(0,1)(1,.5)(1,1.5)
\pspolygon[fillcolor=darkgray](2,2)(2,1)(1,.5)(1,1.5)
\pspolygon[fillcolor=lightgray](-2,2)(-2,1)(-1,.5)(-1,1.5)
\pspolygon[fillcolor=darkgray](0,2)(0,1)(-1,.5)(-1,1.5)
\pspolygon[fillcolor=gray](-1,4.5)(0,4)(1,4.5)(0,5)
\pspolygon[fillcolor=gray](0,3)(1,2.5)(2,3)(1,3.5)
\pspolygon[fillcolor=gray](-2,2)(-1,1.5)(0,2)(-1,2.5)
\psset{fillstyle=none}
\daline(-2,2)(-2,4)(-1,4.5)
\daline(1,4.5)(2,4)(2,3)
\daline(-1,.5)(0,0)(1,.5)(0,1)
\pcline[linestyle=none](-2,1)(-2,4)\Aput{$c$}
\pcline[linestyle=none](0,0)(-2,1)\Aput{$b$}
\pcline[linestyle=none](0,0)(2,1)\Bput{$a$}
\endpspicture}
\rput[b](4.5,0){\pspicture(-2,0)(2,5)
\pspolygon(-1,4.5)(-1,3.5)(0,3)(0,4)
\pspolygon(1,4.5)(1,3.5)(0,3)(0,4)
\pspolygon(-1,3.5)(-1,2.5)(0,2)(0,3)
\pspolygon(0,3)(0,2)(1,1.5)(1,2.5)
\pspolygon(2,3)(2,2)(1,1.5)(1,2.5)
\pspolygon(0,2)(0,1)(1,.5)(1,1.5)
\pspolygon(2,2)(2,1)(1,.5)(1,1.5)
\pspolygon(-2,2)(-2,1)(-1,.5)(-1,1.5)
\pspolygon(0,2)(0,1)(-1,.5)(-1,1.5)
\pspolygon(-1,4.5)(0,4)(1,4.5)(0,5)
\pspolygon(0,3)(1,2.5)(2,3)(1,3.5)
\pspolygon(-2,2)(-1,1.5)(0,2)(-1,2.5)
\psline(-2,2)(-2,4)(-1,4.5)
\psline(1,4.5)(2,4)(2,3)
\psline(-1,.5)(0,0)(1,.5)(0,1)
\psline(-2,3)(-1,3.5)
\pcline[linestyle=none](-2,1)(-2,4)\Aput{$c$}
\pcline[linestyle=none](0,0)(-2,1)\Aput{$b$}
\pcline[linestyle=none](0,0)(2,1)\Bput{$a$}
\endpspicture}
\psline{->}(-1,2.5)(1,2.5)
\endpspicture
\end{fullfigure}

Such a lozenge tiling is equivalent to a matching in a hexagonal mesh
(technically known as chicken wire) graph $Z(a,b,c)$ (\Fig{f:graph}).

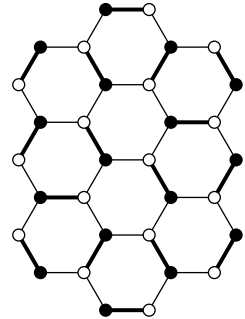
\begin{fullfigure}{f:graph}{The graph $Z(2,2,3)$ and a matching}
\psset{xunit=.289cm,yunit=.5cm}
\pspicture(-6,1)(6,9)
\multips(-4,2)(0,2){4}{\qline(0,0)(2,0)}  
\multips(-1,1)(0,2){5}{\qline(0,0)(2,0)}
\multips( 2,2)(0,2){4}{\qline(0,0)(2,0)}
\multips(-4,2)(0,2){3}{\qline(0,0)(-1,1)} 
\multips(-1,1)(0,2){4}{\qline(0,0)(-1,1)}
\multips( 2,2)(0,2){4}{\qline(0,0)(-1,1)}
\multips( 5,3)(0,2){3}{\qline(0,0)(-1,1)}
\multips( 4,2)(0,2){3}{\qline(0,0)( 1,1)} 
\multips( 1,1)(0,2){4}{\qline(0,0)( 1,1)}
\multips(-2,2)(0,2){4}{\qline(0,0)( 1,1)}
\multips(-5,3)(0,2){3}{\qline(0,0)( 1,1)}
\thline(-1,1)( 1,1) \thline(-5,3)(-4,2)
\thline(-2,2)(-1,3) \thline( 1,3)( 2,2)
\thline( 4,2)( 5,3) \thline(-4,4)(-2,4)
\thline( 1,5)( 2,4) \thline( 4,4)( 5,5)
\thline(-5,5)(-4,6) \thline(-2,6)(-1,5)
\thline( 2,6)( 4,6) \thline(-5,7)(-4,8)
\thline(-2,8)(-1,7) \thline( 1,7)( 2,8)
\thline( 4,8)( 5,7) \thline(-1,9)( 1,9)
\multips(-5,3)(0,2){3}{\whitedot}\multips(-2,2)(0,2){4}{\whitedot}
\multips( 1,1)(0,2){5}{\whitedot}\multips( 4,2)(0,2){4}{\whitedot}
\multips( 5,3)(0,2){3}{\blackdot}\multips( 2,2)(0,2){4}{\blackdot}
\multips(-1,1)(0,2){5}{\blackdot}\multips(-4,2)(0,2){4}{\blackdot}
\endpspicture
\end{fullfigure}

The total number of plane partitions in an $a \times b \times c$ box is given
by MacMahon's formula:
\begin{equation}
N(a,b,c) = \frac{H(a)H(b)H(c)H(a+b+c)}{H(a+b)H(a+c)H(b+c)}, \label{e:macmahon}
\end{equation}
where
$$H(n) = (n-1)!(n-2)!\ldots3!2!1!$$
is the hyperfactorial function. One can enumerate plane partitions with
symmetry, each corresponding to matchings which are invariant under some group
action $G$ on $Z(a,b,c)$.  If $G$ acts freely, these are just the matchings of
the quotient graph $Z(a,b,c)/G$.  When $G$ does not act freely, there is always
a way to modify the quotient graph so that its matchings are equinumerous with
the invariant matchings of $Z(a,b,c)$. The number of plane partitions in each
symmetry class has a formula similar to equation~\eqref{e:macmahon}.

\begin{table}[htb]
\begin{tabular}{lll}
Number & Kind & Acronym \\ \hline
$N_3(a,a,a)$ & Cyclically symmetric & CSPP \\
$N_5(a,b,c)$ & Self-complementary & SCPP \\
$N_9(2a,2a,2a)$ & \begin{tabular}{ll}Cyclically symmetric, \\ self-complementary\end{tabular} & CSSCPP \\
$N_6(2a,b,b)$ & Transpose-complement & TCPP \\ 
$N_8(2a,2a,2a)$ & 
\begin{tabular}{ll}Cyclically symmetric, \\ transpose-complement\end{tabular} & CSTCPP
\end{tabular}
\caption{Plane partitions and their numbers}
\label{t:pp}
\end{table}

A plane partition is {\em cyclically symmetric} if the corresponding matching
is invariant under rotation by $120^\circ$.  It is {\em self-complementary} if
the matching is invariant under rotation by $180^\circ$.  It is {\em
transpose-complement} if the matching is invariant under a color-preserving
reflection. The particular symmetry classes of plane partitions based on these
symmetries that we will consider are given in Table~\ref{t:pp}.  In the table
the parameters of an enumeration $N_i(a,b,c)$ refer to the three dimensions of
the box.  For example $N_6(2a,b,b)$ is the number of TCPPs in a $2a \times b
\times b$ box.

\begin{fullfigure}{f:tcpp}{A deleted quotient graph for TCPPs}
\psset{xunit=.289cm,yunit=.5cm}
\pspicture(-1,0)(17,7)
\psline(-1,0)(-1,4)(8,7)(17,4)(17,0)
\daline(-1,0)(17,0)

\multips( 0,1)(0,2){2}{\thline(0,0)(1,1)}
\multips( 3,2)(0,2){2}{\qline(0,0)(1,1)}
\multips( 6,1)(0,2){3}{\qline(0,0)(1,1)}
\multips( 9,2)(0,2){2}{\qline(0,0)(1,1)}
\multips(12,1)(0,2){2}{\qline(0,0)(1,1)}
\multips(15,2)(0,2){1}{\qline(0,0)(1,1)}

\multips( 1,2)(0,2){1}{\qline(0,0)(-1,1)}
\multips( 4,1)(0,2){2}{\qline(0,0)(-1,1)}
\multips( 7,2)(0,2){2}{\qline(0,0)(-1,1)}
\multips(10,1)(0,2){3}{\qline(0,0)(-1,1)}
\multips(13,2)(0,2){2}{\qline(0,0)(-1,1)}
\multips(16,1)(0,2){2}{\thline(0,0)(-1,1)}

\multips(0,0)(6,0){3}{\daline(0,1)(1,0)(3,0)(4,1)
    \pscircle[linestyle=dashed,fillstyle=solid,fillcolor=lightgray](1,0){2.5pt}
    \pscircle[linestyle=dashed,fillstyle=solid,fillcolor=white](3,0){2.5pt}}

\multips( 1,2)(0,2){2}{\qline(0,0)(2,0)\blackdisk(0,0){2.5pt}\whitedisk(2,0){2.5pt}}
\multips( 4,1)(0,2){3}{\qline(0,0)(2,0)\blackdisk(0,0){2.5pt}\whitedisk(2,0){2.5pt}}
\multips( 7,2)(0,2){3}{\qline(0,0)(2,0)\blackdisk(0,0){2.5pt}\whitedisk(2,0){2.5pt}}
\multips(10,1)(0,2){3}{\qline(0,0)(2,0)\blackdisk(0,0){2.5pt}\whitedisk(2,0){2.5pt}}
\multips(13,2)(0,2){2}{\qline(0,0)(2,0)\blackdisk(0,0){2.5pt}\whitedisk(2,0){2.5pt}}
\whitedisk( 0,1){2.5pt}\whitedisk( 0,3){2.5pt}
\blackdisk(16,1){2.5pt}\blackdisk(16,3){2.5pt}
\endpspicture
\end{fullfigure}

\begin{fullfigure}{f:cstcpp}{A deleted quotient graph for CSTCPPs}
\psset{xunit=.5cm,yunit=.289cm}
\pspicture(0,-9)(12,9)

\daline(9,-9)(0,0)(9,9)
\psline(9,-9)(12,0)(9,9)

\multips(0,0)(3,-3){3}{\daline(1, 1)(1,-1)(2,-2)(3,-1)
    \pscircle[linestyle=dashed,fillstyle=solid,fillcolor=lightgray](1,-1){2.5pt}
    \pscircle[linestyle=dashed,fillstyle=solid,fillcolor=white](2,-2){2.5pt}}

\multips(0,0)(3, 3){3}{\daline(1,-1)(1, 1)(2, 2)(3, 1)
    \pscircle[linestyle=dashed,fillstyle=solid,fillcolor=lightgray](2, 2){2.5pt}
    \pscircle[linestyle=dashed,fillstyle=solid,fillcolor=white](1, 1){2.5pt}}

\multips(6,-4)(2,0){2}{\psline(0,0)(1,-1)(2,0)}
\multips(3,-1)(2,0){4}{\psline(0,0)(1,-1)(2,0)}
\multips(3, 1)(2,0){4}{\psline(0,0)(1, 1)(2,0)}
\multips(6, 4)(2,0){2}{\psline(0,0)(1, 1)(2,0)}
\thline( 9,-7)( 9,-5) \thline(10,-4)(10,-2) \thline(11,-1)(11, 1)
\thline(10, 2)(10, 4) \thline( 9, 5)( 9, 7)

\multips(9,-7)(2,0){1}{\qline(0,0)(0,2)\blackdisk(0,0){2.5pt}\whitedisk(0,2){2.5pt}}
\multips(6,-4)(2,0){3}{\qline(0,0)(0,2)\blackdisk(0,0){2.5pt}\whitedisk(0,2){2.5pt}}
\multips(3,-1)(2,0){5}{\qline(0,0)(0,2)\blackdisk(0,0){2.5pt}\whitedisk(0,2){2.5pt}}
\multips(6, 2)(2,0){3}{\qline(0,0)(0,2)\blackdisk(0,0){2.5pt}\whitedisk(0,2){2.5pt}}
\multips(9, 5)(2,0){1}{\qline(0,0)(0,2)\blackdisk(0,0){2.5pt}\whitedisk(0,2){2.5pt}}
\whitedisk(7,-5){2.5pt}
\whitedisk(4,-2){2.5pt}
\blackdisk(4, 2){2.5pt}
\blackdisk(7, 5){2.5pt}
\endpspicture
\end{fullfigure}

The invariant matchings corresponding to the first three symmetry classes are
equivalent to matchings of $Z(a,b,c)/G$ for suitable $G$.  The partitions in
the last two symmetry classes correspond to matchings of $Z(a,b,c)//G$, a graph
which is obtained from the quotient graph by deleting those vertices that have
a non-trivial stabilizer in $G$ (Figures~\ref{f:tcpp} and \ref{f:cstcpp}). In
each case let $Z_i(a,b,c)$ be the graph corresponding to $N_i(a,b,c)$.   In the
last two cases, some of the vertices are {\em vestigial}, meaning that they can
only be matched in one way.  (These vertices are matched with bold edges in the
figures.)  Let $Z'_i(a,b,c)$ be the graph $Z(a,b,c)$ with vestigial vertices
removed.

The $q$-enumeration $N(a,b,c)_q$ of plane partitions is given by the natural
$q$-analogue of equation~\eqref{e:macmahon}, namely the one where each
factorial is replaced by a $q$-factorial.  Here the $q$-weight of a plane partition is $q^k$ if there
are $k$ cubes.  The $q$-enumeration is also the determinant of $Z(a,b,c)$ if it
is weighted with curvature $q$ in each hexagonal face. Three of the other
symmetry classes can also be $q$-enumerated, but we will consider only
the $q$-enumeration of cyclically symmetric plane partitions, where
as with unrestricted plane partitions, the weight of each cube is $q$.

The six symmetry classes that do involve complementation have no obvious
$q$-enumerations, but they do have natural $-1$ enumerations.  In these cases
the symmetry group $G$ acts on individual cubes just as it does without
complementation.  But whereas in the four classes without complementation each
orbit of cubes is either filled or empty, in the six classes with
complementation each orbit is always half-filled.  Nonetheless, there is a
natural move between symmetric plane partitions which replaces half of an orbit
of cubes by the opposite half.  Any two plane partitions in a given symmetry
class differ by either an odd or an even number of moves.  This defines a
relative sign between them.  The corresponding signed enumeration is a natural
generalization of $q$-enumeration with $q=-1$. We conjecture that there is a
product formula for the $-1$-enumeration of every symmetry class of plane
partitions.

Stembridge \cite{Stembridge:strange} has found interesting patterns in signed
enumerations of cyclically symmetric plane partitions called {\em strange}
enumerations.  In a strange enumeration, the weight of each plane partition is
the product of the weights of cubes. Each cube has weight $1$ or $-1$ depending
on its position in the box.  Consider plane partitions in the box $[0,a]^3$ in
Cartesian coordinates, and denote each unit cube by the coordinates of the
corner farthest away from the origin.  We will consider strange enumerations
$N_3(a,a,a)_{s,t,u}$, where $s$ is the sign (or weight) of cubes $(i,i,i)$, $t$
is the sign of cubes $(i,j,j)$ with $i \ne j$, and $u$ is the sign of cubes
$(i,j,k)$ with $i$, $j$, and $k$ all different. Similarly, let $N(2a,2b,2c)_{s,t}$ be a
signed (or weighted) enumeration of unrestricted plane partitions in a $2a
\times 2b \times 2c$ box, where the cubes $(a+i,b+i,c+i)$ (where $i$ can be
positive or negative) have weight $s$ and the other cubes have weight $t$.

\begin{fullfigure}{f:penrose}{A Penrose-style lozenge tiling}
\psset{xunit=.5cm,yunit=.2887}
\pspicture(0,0)(5,10)
\pspolygon[fillstyle=solid,fillcolor=lightgray](2,4)(3,5)(2,6)
\qline(0,2)(0,8)
\qline(1,1)(1,7)
\qline(2,0)(2,2)\qline(2,4)(2,8)
\qline(3,1)(3,3)\qline(3,5)(3,7)
\qline(4,2)(4,6)
\qline(5,3)(5,7)

\qline(0,2)(2,0)
\qline(0,4)(2,2)
\qline(0,6)(1,5)\qline(2,4)(3,3)
\qline(0,8)(1,7)\qline(2,6)(4,4)
\qline(1,9)(4,6)
\qline(2,10)(5,7)

\qline(2,0)(5,3)
\qline(2,2)(5,5)
\qline(1,3)(3,5)\qline(4,6)(5,7)
\qline(1,5)(2,6)\qline(3,7)(4,8)
\qline(1,7)(3,9)
\qline(0,8)(2,10)
\endpspicture
\end{fullfigure}

Finally, let $a$, $b$, and $c$ have the same parity.  We will consider lozenge
tilings of an $a$, $b+1$, $c$, $a+1$, $b$, $c+1$ hexagon with the middle unit
triangle removed.  Such tilings resemble a Penrose impossible triangle
(\Fig{f:penrose}).  Let $N_P(a,b,c)$ be the number of these and let $Z_P(a,b,c)$
be the corresponding graph. There is a product formula for $N_P(a,b,c)$
\cite{Propp:twenty} which has been proved by Okada and Krattenthaler
\cite{OK:punctured}. A very short argument has been found by Ciucu for the case
$b=c$ \cite{Ciucu:punctured} and generalized \cite{Ciucu:macmahon}.

\subsection{Various relations}
\label{s:various}

The properties of the permanent-determinant method yield various relations
between the enumerations of plane partitions defined above.  Many enumerations
of plane partitions are {\em round}. An integer is round if it is a product of
relatively small numbers.  A polynomial is round if it is a product of
cyclotomic polynomials of low degree.  A round enumeration in combinatorics
almost always has an explicit product formula.

Rotational symmetry of $Z(a,a,a)$ implies that
$$N(a,a,a)_{q^3} = N_3(a,a,a)_q N_3(a,a,a)_{\omega q}
N_3(a,a,a)_{\omega^2 q},$$
where $\omega$ is a cube root of unity.  This relation together with MacMahon's
formula for $N(a,a,a)_q$ implies that $N_3(a,a,a)_q$ is round, but it does not
imply the explicit formula for $N_3(a,a,a)_q$ conjectured by MacDonald and
proved by Mills, Robbins, and Rumsey \cite{MRR:macdonald}.

Rotational symmetry also tells us that
\begin{align*}
N(a,a,a)_{1,-1} =\ & N_3(a,a,a)_{1,-1,-1} \\
& N_3(a,a,a)_{\omega,-1,-1} \\
& N_3(a,a,a)_{\omega^2,-1,-1}.
\end{align*}
The factor $N_3(a,a,a)_{1,-1,-1}$ is a strange enumeration.  This case is
entirely analogous to $N_3(a,a,a)_q$ because we know that $N(a,a,a)_{1,-1}$ is
round (see below).  Therefore $N_3(a,a,a)_{1,-1,-1}$ is round also, but the
explicit formula conjectured by Stembridge is still open.

Reflection symmetry in $Z(2a,b,b)$ tells us that
$$N(2a,b,b) = N_6(2a,b,b)N_6(2a,b,b)_2,$$
and reflection symmetry in $Z(2a+1,b,b)$ tells us that
$$N(2a+1,b,b) = 2N_6(2a,b,b)N_6(2a,b,b)_2.$$
Here the second factor is a certain 2-enumeration of TCPP's (a weighting where
certain faces have curvature 2). These two relations imply the known formula
for $N_6(2a,b,b)$ \cite{Stanley:symmetries}.

The method of Section~\ref{s:reverse} establishes the identities
\begin{align*}
|N_1(2a,2b,2c)_{-1,1}| & = N_5(2a,2b,2c)^2 \\
|N_3(2a,2a,2a)_{-1,1,1}| & = N_9(2a,2a,2a)^2.
\end{align*}
The second identity is one of Stembridge's strange enumerations. Two other
strange enumerations,$|N_3(2a,2a,2a)_{-1,-1,1}|$ and
$|N_3(2a,2a,2a)_{-1,1,-1}|$, can be recognized as perfect squares by the same
method.  However, these enumerations are not round numbers; in general they
have large prime factors.

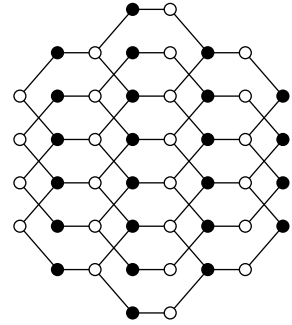
\begin{fullfigure}{f:overlay}{Two copies of $Z(2,2,2)$ overlaid}
\psset{xunit=.333cm,yunit=.577cm}
\pspicture(-6,1)(6,8)
\multips(-4,2)(0,1){6}{\qline(.25,0)(1.75,0)}  
\multips(-1,1)(0,1){8}{\qline(.25,0)(1.75,0)}
\multips( 2,2)(0,1){6}{\qline(.25,0)(1.75,0)}
\multips(-4,2)(0,1){4}{\qline(.25,0)(-1.25,1)} 
\multips(-1,1)(0,1){6}{\qline(.25,0)(-1.25,1)}
\multips( 2,2)(0,1){6}{\qline(.25,0)(-1.25,1)}
\multips( 5,3)(0,1){4}{\qline(.25,0)(-1.25,1)}
\multips( 4,2)(0,1){4}{\qline(-.25,0)(1.25,1)} 
\multips( 1,1)(0,1){6}{\qline(-.25,0)(1.25,1)}
\multips(-2,2)(0,1){6}{\qline(-.25,0)(1.25,1)}
\multips(-5,3)(0,1){4}{\qline(-.25,0)(1.25,1)}
\multips(-5.25,3)(0,1){4}{\whitedot}\multips(-2.25,2)(0,1){6}{\whitedot}
\multips( 0.75,1)(0,1){8}{\whitedot}\multips( 3.75,2)(0,1){6}{\whitedot}
\multips( 5.25,3)(0,1){4}{\blackdot}\multips( 2.25,2)(0,1){6}{\blackdot}
\multips(-0.75,1)(0,1){8}{\blackdot}\multips(-3.75,2)(0,1){6}{\blackdot}
\endpspicture
\end{fullfigure}

The construction in Section~\ref{s:force} applies to plane partitions as
follows:  Let $Z$ be the disjoint union of two copies of the plane partition
graph $Z(a,b,c)$ with flat weighting.  Arrange $Z$ so that most of the
horizontal edges of one copy are concentric with hexagons of the other copy, is
in \Fig{f:overlay}. Now replace all of the crossings by butterflies.  Most of
the long horizontal edges of the butterflies cancel in pairs (they are multiple
edges with opposite weight).  The new graph, after a little vertex splitting,
is exactly the graph $Z(2a,2b,2c)$, but the weighting is no longer flat. Rather
it has curvature $-1$ in every face.  Since
$$\det(Z) = N(a,b,c)^2,$$
the conclusion is that
$$N(2a,2b,2c)_{-1} = N(a,b,c)^2.$$
This relation is independently a corollary of MacMahon's $q$-enumeration of plane
partitions.  Via Stembridge's $q=-1$ phenomenon, it is also related to the
enumeration of self-complementary plane partitions. Since the original graph
$Z$ could have been given a weighting with indeterminate curvature instead of a
flat one, it generalizes to a new relation between generating functions for
plane partitions in a $2a \times 2b \times 2c$ box and plane partitions in an
$a \times b \times c$ box.

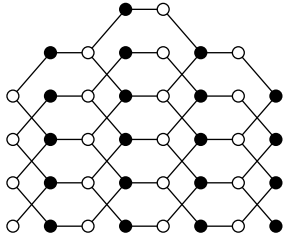
\begin{fullfigure}{f:overlay6}{Overlaying $Z_6(4,2,2)$ and $Z'_6(4,3,3)$}
\psset{xunit=.333cm,yunit=.577cm}
\pspicture(-5,0)(5,5)
\multips(-4,0)(0,1){5}{\qline(.25,0)(1.75,0)}  
\multips(-1,0)(0,1){6}{\qline(.25,0)(1.75,0)}
\multips( 2,0)(0,1){5}{\qline(.25,0)(1.75,0)}
\multips(-4,0)(0,1){3}{\qline(.25,0)(-1.25,1)} 
\multips(-1,0)(0,1){4}{\qline(.25,0)(-1.25,1)}
\multips( 2,0)(0,1){5}{\qline(.25,0)(-1.25,1)}
\multips( 5,0)(0,1){4}{\qline(.25,0)(-1.25,1)}
\multips( 4,0)(0,1){3}{\qline(-.25,0)(1.25,1)} 
\multips( 1,0)(0,1){4}{\qline(-.25,0)(1.25,1)}
\multips(-2,0)(0,1){5}{\qline(-.25,0)(1.25,1)}
\multips(-5,0)(0,1){4}{\qline(-.25,0)(1.25,1)}
\multips(-5.25,0)(0,1){4}{\whitedot}\multips(-2.25,0)(0,1){5}{\whitedot}
\multips( 0.75,0)(0,1){6}{\whitedot}\multips( 3.75,0)(0,1){5}{\whitedot}
\multips( 5.25,0)(0,1){4}{\blackdot}\multips( 2.25,0)(0,1){5}{\blackdot}
\multips(-0.75,0)(0,1){6}{\blackdot}\multips(-3.75,0)(0,1){5}{\blackdot}
\endpspicture
\end{fullfigure}

Forcing planarity also establishes the identities
\begin{align*}
N(2a+1,2b,2c)_{-1} =\ & N(a,b,c)N(a+1,b,c) \\
N(2a+1,2b+1,2c)_{-1} =\ & N(a,b+1,c)N(a+1,b,c) \\
N_6(4a,2b,2b)_{-1} =\ & N_6(2a,b,b)^2 \\
N_6(4a+2,2b+1,2b+1)_{-1} =\ & N_6(2a,b+1,b+1)\\ & N_6(2a+2,b,b) \\
N_6(4a,2b+1,2b+1)_{-1} =\ & N_6(2a,b,b) \\ & N_6(2a,b+1,b+1).
\end{align*}
The last three identities require some manipulations with
vestigial vertices.  For example, the last identity
is derived by overlaying $Z_6(2a,b,b)$ and $Z'_6(2a,b+1,b+1)$
and then forcing planarity (\Fig{f:overlay6}).

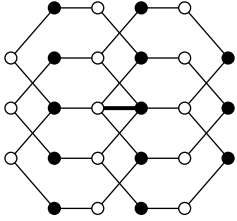
\begin{fullfigure}{f:overlayp}{Two joined copies of $N_P(1,1,1)$}
\psset{xunit=.385cm,yunit=.667cm}
\pspicture(0,0)(7,4)
\multips(1,0)(0,1){5}{\qline(.25,0)(1.75,0)}
\multips(4,0)(0,1){5}{\qline(.25,0)(1.75,0)}
\multips(1,0)(0,1){3}{\qline(.25,0)(-1.25,1)}
\multips(4,0)(0,1){4}{\qline(.25,0)(-1.25,1)}
\multips(7,1)(0,1){3}{\qline(.25,0)(-1.25,1)}
\multips(0,1)(0,1){3}{\qline(-.25,0)(1.25,1)}
\multips(3,0)(0,1){4}{\qline(-.25,0)(1.25,1)}
\multips(6,0)(0,1){3}{\qline(-.25,0)(1.25,1)}
\thline(2.75,2)(4.25,2)
\multips(-.25,1)(0,1){3}{\whitedot}
\multips(2.75,0)(0,1){5}{\whitedot}
\multips(5.75,0)(0,1){5}{\whitedot}
\multips(1.25,0)(0,1){5}{\blackdot}
\multips(4.25,0)(0,1){5}{\blackdot}
\multips(7.25,1)(0,1){3}{\blackdot}
\endpspicture
\end{fullfigure}

A more complicated use of \Sec{s:force} establishes the identities
\begin{equation}
|N(2a,2b,2c)_{1,-1}| = N(a,b,c)N(a,b,c)_{-1}^2 \label{e:combo}
\end{equation}
and
\begin{equation}
|N(2a+1,2b+1,2c+1)_{1,-1}| = 2N_P(a,b,c)^2. \label{e:penrose}
\end{equation}
To derive equation~\eqref{e:combo}, take two copies of $Z(a,b,c)$, one flat and
one with curvature $-1$ at the center face, and overlay them as in
\Fig{f:overlay}.  The first copy has determinant $N(a,b,c)$, while by
\Sec{s:reverse}, the second copy has determinant $N(a,b,c)_{-1}^2$.  For
equation~\eqref{e:penrose}, take two copies of $Z_P(a,b,c)$ overlaid, but
instead of deleting the middle vertices, connect them by an edge of weight 2,
as in \Fig{f:overlayp}.  (In the figure, the connecting edge of weight 2 is in
bold.) If each copy of $Z_P(a,b,c)$ is weighted so that it would be flat if its
middle vertex were removed, then the determinant of the combined graph is
$2N_P(a,b,c)^2$, since every matching has to contain the middle edge.  But
forcing planarity converts this graph to $Z(2a+1,2b+1,2c+1)$ with a certain
weighting.  The determinant of this graph is the left side of
equation~\eqref{e:penrose}, as desired.

Symmetry of $Z(a,a,a)$ yields the factorization:
\begin{align*}
N(a,a,a)_{1,-1} =\ &N_3(a,a,a)_{1,-1,-1} \\
& N_3(a,a,a)_{\omega,-1,-1} \\
& N_3(a,a,a)_{\omega^2,-1,-1}.
\end{align*}
This implies that $N_3(a,a,a)_{1,-1,-1}$ is a round number,
but does not establish its explicit formula.

Finally, plane partitions in an $a \times b \times c$ box are a good
example of the connection between the permanent-determinant method and the
Gessel-Viennot method.  Here there are three Gessel-Viennot matrices, of order
$a \times a$, $b \times b$, and $c \times c$.  They are also known as
Carlitz matrices.

\begin{question}  What is the cokernel of a Kasteleyn matrix for plane
partitions in an $a \times b \times c$ box?
\end{question}

By Gessel-Viennot, if
$$\min(a,b,c) = 1,$$
the cokernel is cyclic.  But if
$$a=b=c=2,$$
the cokernel is the non-cyclic abelian group $\Z/2 \times \Z/10$.

\section{Historical note}
\label{s:historical}

Kasteleyn \cite{Kasteleyn:crystal} gives an excellent exposition of both the
Hafnian-Pfaffian method and its application to problems in statistical 
mechanics.  His bibliography constitutes a nearly complete chronicle of the
series of papers that led to its discovery.  In this section, we give a 
summary of that chronicle.

Onsager's solution of the Ising model \cite{Onsager:crystal} was the first 
convincing analysis of any interesting statistical mechanical model of a
crystal.  However, this paper contained difficult mathematics which left open
the possibility of a simpler approach.   With this motivation, Kac and Ward
\cite{KW:ising} announced a new solution eight years later, in which the
partition function of the Ising model was expressed as the determinant of a
relatively small matrix.  The permutation expansion of the determinant was
interpreted as a weighted enumeration.  This paper attracted the attention of
many strong physicists and mathematicians \cite{NM:ising}, but it contained a
serious combinatorial error.  Feynman recast the calculations in a correct but
incomplete form in unpublished work.  He replaced the error by a conjecture
which was popularized by Harary, who was organizing a book on applications of
graph theory to physics. Sherman proved the conjecture \cite{Sherman:ising} and
generalized both it and the Kac-Ward determinant to the Ising model for
arbitrary plane graphs.  Harary eventually substituted his reference to
Sherman's method by Kasteleyn's survey \cite{Kasteleyn:crystal}.

Independently of Sherman, Hurst and Green \cite{HG:ising} found a cleaner
approach to the square-lattice Ising model in which the partition function is
the Pfaffian of the antisymmetric incidence matrix of a nearly planar,
weighted graph. The paper of Hurst and Green inspired Fisher and
Temperley, and independently Kasteleyn, to enumerate perfect matchings in a
square lattice with a special case of the Hafnian-Pfaffian method.  Fisher saw
that the method also works for finite hexagonal meshes (which already includes
some enumerations of current interest).  Independently Kasteleyn
generalized the method to all plane graphs.

Later Percus \cite{Percus:dimer} saw that since the square lattice is
bipartite, the Hafnian-Pfaffian method reduces to the permanent-determinant
method.  With hindsight it seems that the only reason that Pfaffians and not
determinants were used for the square lattice dimer model  was that they arise
naturally in the planar Ising model.  (Mathematical physicists call weighted
enumerations of perfect matchings {\em dimer models}.)  Much later still,
Little showed that $K_{3,3}$-free graphs are Pfaffian \cite{Little:k33free}.  

These results seemed to complete the story of the permanent-determinant method,
until interest was renewed more recently in the context of enumerative
combinatorics
\cite{Kuperberg:eklp1,Kuperberg:eklp2,Jockusch:perfect,Kuperberg:perdet,%
Kuperberg:oneroof,Propp:twenty,Yang:thesis,Tesler:non-oriented}.  In the author's opinion, one of
the nicest new observations about the method is the fact that some of the
minors of the inverse of a large Kasteleyn matrix are the probabilities of
local configurations of edges \cite{Kenyon:local}. For this reason and others,
the method could well have some bearing on arctic circle phenomena in domino
and lozenge tilings \cite{JPS:arctic,CLP:typical}. Another recent development
is a satisfactory classification of bipartite Pfaffian graphs
\cite{RST:Pfaffian}.

The Gessel-Viennot method was discovered by Gessel and Viennot
\cite{GV:binomial,GV:partitions} in the context of enumerative combinatorics,
and completely independently of the permanent-determinant method.   The method
was independently anticipated by Lindstr\"om \cite{Lindstrom:vector} in the
context of matroid theory.  It has been widely influential in enumeration of
plane partitions, even more so than the permanent-determinant method
\cite{Andrews:tsscpp,AB:determinant,Ciucu:macmahon,CK:plane,%
Krattenthaler:determinant,Okada:generating,OK:punctured,Stembridge:enumeration}.
Yet another determinant variation which is likely related to the
permanent-determinant method was recently found by Helfgott and Gessel
\cite{HG:diamonds}.


\end{document}